\documentstyle[12pt,epsf,titlepage]{article}
\newcommand{\del}{\triangle}
\newcommand{\pd}{\partial}
\newcommand{\ol}{\overline}
\newcommand{\wdtld}{\widetilde}
\newcommand{\xa}{x_\alpha}
\newcommand{\xb}{x_\beta}
\newcommand{\ua}{u_\alpha}
\newcommand{\ub}{u_\beta}
\newcommand{\fua}{\overline{u}_\alpha}

\newcommand{\bx}{{\bf x}}
\newcommand{\bs}{{\bf s}}

\renewcommand{\baselinestretch}{2}
\begin{document}
\begin{titlepage}
\setlength{\topmargin}{0in}
\begin{center}
\thispagestyle{plain}
\pagenumbering{roman}
{\LARGE A Modified Smagorinsky Subgrid Scale Model for the Large Eddy Simulation
of Turbulent Flow}
\\[0.25in]
BY
\\
{\large Tommy Kunhung Kim}
\\

B.A.~(California State University, Fullerton) 1990
\\ \vspace{-0.25in}
M.A.~(California State University, Fullerton) 1992
\\
{\large DISSERTATION}
\\
Submitted in partial satisfaction of the requirements for the degree of
\\
DOCTOR OF PHILOSOPHY
\\
in
\\
Applied Mathematics
\\
in the 
\\
OFFICE OF GRADUATE STUDIES
\\
of the
\\
UNIVERISTY OF CALIFORNIA
\\
DAVIS
\\

Approved: \hrulefill \hspace{1.25in}~
\\[.1in]
\hspace{.79in}\hrulefill \hspace{1.25in}~ 
\\[.1in]
\hspace{.79in}\hrulefill \hspace{1.25in}~
\\

Committee in Charge
\\
2001
\end{center}
\end{titlepage}

\begin{abstract}
\thispagestyle{plain}
\pagenumbering{roman}
\setcounter{page}{2}
In the field of Large Eddy Simulation, the Smagorinsky subgrid scale model
(in some form) is the most commonly accepted and used subgrid scale model. The purpose of
this paper is to address the main weakness of the Smagorinsky model, its
poor performance near the wall. The goal is to establish a model that
corrects the Smagorinsky model near the walls while at the same time minimizing
the computational overhead. A version the Dynamic Subgrid Scale model is also
incorporated into the finite element code to facilitate comparisons with the
new model near the walls. 

One of the unique characteristics of Large Eddy Simulations as compared to 
other methods of dealing with turbulent flows is the idea of filtering. In
this paper we define what a filter is and also address an issue related to
filters; the error that results when the filtering and differential
operations are interchanged. This error is studied under the context of the 
Finite Element Method which allows us to focus on the function being filtered
rather than the filter kernal function, which has been the usual approach in
studying this error.
\end{abstract}

\thispagestyle{plain}
\pagenumbering{roman}
\setcounter{page}{3}
\tableofcontents
\thispagestyle{plain}
\newpage

\section{Introduction}
\pagenumbering{arabic}

In the study of turbulence, three distinct computational approaches 
have been developed to simulate turbulent flows. Direct numerical simulation 
(DNS), as the name implies, is the most straight forward approach to the 
simulation of turbulent flows. The DNS approach is to simply solve the 
Navier-Stokes equations using some numerical method. The problem with this 
approach is that to 
succesfully solve the Navier-Stokes equations, the computational domain must
be large enough to contain the largest scales of motion, $L$, while the grid 
resolution must also be fine enough to resolve the smallest scales, which are
of the order of the Kolmogorov microscale, $\eta = (\nu^3/\epsilon)^{1/4}$,
where $\eta$ is the dissipation rate per unit mass and $\nu$ is the 
kinematic viscosity. So, the number of grid points necessary in each direction
can be estimated as $\frac{L}{\eta}$ and using the well established 
relationship between $\frac{L}{\eta}$ and the Reynolds number, $Re$,
(see Tennekes and Lumely, \cite{lumley}) we have that
\[ N = \left( \frac{L}{\eta} \right)^3 \sim Re^{\frac{9}{4}}, \]
where $N$ is the number of grid points. Hence, it is clear that DNS is only
practical for low Reynolds number flows. At the other end of the spectrum
from the DNS approach is the Reynolds averaged simulations (RAS). In RAS,
the flow (velocity field and pressure for incompressible flows) is broken up
into a statistically steady (mean) portion and fluctuations. The 
mean flow is then solved for, while the effects of the fluctuations on the flow
are modelled. By modeling the fluctuation effects, the RAS approach models 
all the turbulent interactions whereas the DNS approach completely resolves 
them. One of the main stumbling blocks for RAS is that no one model has been
found for the fluctuation effects that can be used for all turbulent flows.
It seems unlikely that such a model would exist given the complexity of 
turbulent flows and the fact that turbulence is a property of the flow rather
than the fluid. The third approach to the computation of turbulent flows, 
known as Large-Eddy simulation (LES), is often viewed as the intermediary 
approach between DNS and RAS. The motivation behind LES is that since the 
large energy carrying eddies are highly influenced by the boundary conditions,
it should be computationally resolved, while the small eddies or unresolved
scales are modelled. It is hoped that since the small eddies are more 
homogeneous and isotropic, a simpler and more universal model can be used
for the unresolved scales as compared to the RAS models. 

The objective of this study is to develope a new model for the unresolved
scales that overcomes some of the problems associated with the established
models, which are discussed in the coming sections. This study is based on
and, in many ways, an extension of the work begun by Rose McCallen
\cite{mccallen} in which a LES method is incorporated into a finite element
code. In her study, a classical model for the unresolved scales developed
by Smagorinsky \cite{smag} was used to simulate the channel and 
backward facing step flows in two dimensions. In the present study, this
model plus a version of a newer model first developed by Germano, Piomelli,
Moin, and Cabbot \cite{germano1} and our own model base on the Smagorinsky
model are used to simulate the channel and backward facing step flows in
three dimensions. 

\newpage

\section{LES}

\subsection{Filtering}

In LES, the underlying principle is to seperate the large eddies, which
are computationally resolved, from the small, which are modelled. This 
partitioning of the scales is accomplished by filtering, in the case
of incompressible flow, the velocity field, {\bf u}, and pressure, $P$.
Following the work of Leonard \cite{leonard}, if $f({\bf x})$ is a function
that contains all the scales, then the filter of $f$, denoted as $\ol{f}$, is
defined as the convolution of $f$ with a filter function $G({\bf x})$.
\begin{equation}
\label{defn-filter}
  \ol{f}({\bf x}) = \int_\Gamma G({\bf x-s}) f({\bf s}) \,d{\bf s}
\end{equation}
where $\Gamma$ is the flow volume, $G$ is normalized so that 
\[ \int_\Gamma G({\bf s}) \,d{\bf s} = 1, \] 
and $G(-x) = G(x)$.
This is necessary to insure that when $f$ is constant that $\ol{f} = f$.
Some common filters are the following:
\begin{itemize}
  \item Top--Hat filter
    \[ G({\bf x - s}) = \left\{ \begin{array}{c}
           \frac{1}{\del_f^3} \quad \mbox{for} |{\bf x - s}| < \frac{\del_f}{2}
           \\
             0 \quad \mbox{otherwise}
          \end{array} \right.
     \]
    where $\del_f$ is the filter size. If {\bf u} is the velocity field, then
    the filter of {\bf u}, using the Top-Hat filter, is
    \[ {\bf \ol{u}(x)} = \frac{1}{\del^3_f} 
                         \int^{\frac{\del_f}{2}}_{-\frac{\del_f}{2}}
                              {\bf u(x + s)} \,d{\bf s}
    \]
    which is equivalent to cell volume averaging of Deardorff \cite{deardorff}.
    Now, the Fourier transform of ${\bf \ol{u}(x)}$ is
    \begin{eqnarray}
    \label{spectral-eqn}
      {\bf \hat{\ol{u}}(k)} = \left[ \prod^3_{i=1}
           \frac{\sin\left(k_i \frac{\del_f}{2}\right)}{k_i \frac{\del_f}{2}}
           \right] {\bf \hat{u}(k)},
    \end{eqnarray}
     where ${\bf \hat{u}(k)}$ is the Fourier transform of {\bf u}
     (see Kwak $et~al.$~\cite{kwak}.) So, we
     see that the spectrum of the filtered field contains components of all
     wave numbers. This implies that the Top-Hat filter will not filter out
     all the small scales of the flow. Also, the wave numbers for which 
     the coefficient
     of ${\bf \hat{u}(k)}$ in Eqn.~(\ref{spectral-eqn}) is zero lead to the 
     inverse
     transform being singular. This means that it is not possible to obtain
     the actual spectrum ${\bf \hat{u}(k)}$ from the filtered one, 
     ${\bf \hat{\ol{u}}(k)}$. So the Top-Hat filter is not appropriate when
     spectral results are sought.

  \item Gaussian filter
    \[ G({\bf x-s}) = \left(\frac{6}{\pi \del^2} \right)^{\frac{3}{2}}
                      \exp\left(\frac{-6 ({\bf x - s})^2}{\del_f^2} \right)
    \]
    According to Rogallo and Moin \cite{rogallo}, the Gaussian filter is a 
    good filter choice for filtering in the homogeneous directions
    because it provides a smooth transition between the resolved and 
    subgrid (unresolved) scales and is positive definite in both the physcial
    and wave space. Kwak $et~al.$~\cite{kwak} showed that the Gaussian filter
    results in a filtered field that captures most of the large scale motion
    while most of the small scale motions are removed. This ability of the
    Gaussain filter to seperate the large and small scales, gives it a
    desirable property for a filter in LES. The one obvious flaw or weakness
    of this filter is that it does not have a compact domain. Hence, its use
    in non-homogeneous, wall bounded flows is questionable at best. 

  \item Sharp--Cutoff filter
    \[ G(x) = \prod^3_{i=1} \frac{\sin \left( \frac{\pi(x_i-s_i)}{\del_f}
                                      \right)}
             {\pi (x_i - s_i)}
    \]
    Just as the Top-Hat filter has compact support in grid coordinates, the
    Sharp-Cutoff filter has compact support in the wave space. Also, just
    as the Top-Hat filter looses it compact support when transformed into
    the wave space, the Sharp-Cutoff filter looses its compact support when
    going from the wave space to grid coordinates. So, like the Gaussian filter,
    it is not appropriate for non-homogeneous, wall bounded flows.

\end{itemize}

Although Leonard's definition for filtering is generally accepted, two
distinct approaches to the application of the filter have been developed.
These are mainly influenced by the flow under study and the numerical 
technique used to solve the flow problem. In the case of homogenous turbulence,
the use of spectral methods facilitates the use of an explicit filter. However,
in using finite difference methods, usually for non--homogeneous flows, the
filter is assumed to be implicitly defined by the discretization. An example
of this implicit filter is given by Rogallo and Moin, \cite{rogallo}, where
they note that the second--order central--difference of a function $f$ is the
derivative of the filter of $f$:
\[ \frac{f(x+\del) - f(x-\del)}{2\del} 
  = \frac{d}{dx} \left[\frac{1}{2\del}\int^{x+\del}_{x-\del} f(s) \,ds \right]
  = \frac{d \ol{f}}{dx}, \]
where the filter used is the Top--Hat filter. Some of the reasons for using
an implicit filter rather than explicitly applying a filter are the following:
the difficulty of actually applying a filter; the reduction in grid resolution
resulting from the fact that in finite differencing, the function is only
defined at the grid nodal points and so any filtering would require at least
two nodal points (i.~e.~just like the pressure, the differencing would be 
done using the grid nodes, but the filtered value would be at the cell center);
and finally, finite difference schemes are usually used for non--homogenous 
wall bounded flows which usually have anisotrophic grids (since the filter
width will most likely depend on the grid, this means that the filter width
would be a function of location unlike the homogeneous case where the filter
width remains constant. This leads to further difficulties, as will be
discussed below.)

In applying a filter, Eqn.~(\ref{defn-filter}), to get the 
LES equations, a question that arises 
is whether or not the derivative operator and filtering
are interchangeable, i.~e.~is the derivative of the filter equal to the
filter of the derivative. The answer will depend on the flow under study. In
the case of homogeneous turbulence, where the scales of the large and small
eddies do not depend on the location in space, there would be no reason to use
anything but a constant filter width, $\del$. So, assuming that
$G \to 0$ as $t \to \infty$ and using the 1--dimensional case for notational
simplicity (with the extension into higher dimensions following in the usual
way),
\begin{eqnarray*}
  \frac{\pd \ol{f}(x)}{\pd x} 
  &=& \frac{\pd}{\pd x}
      \int_{-\infty}^\infty f(s) G(x-s)\,ds \\
  &=& \int_{-\infty}^\infty f(s) \frac{\pd}{\pd x} G(x-s)\,ds \\
  &=& \int_{-\infty}^\infty \frac{\pd}{\pd s} f(s) G(x-s)\,ds \\
  &=& \ol{\frac{\pd f}{\pd x}},
\end{eqnarray*}
where integration by parts and the fact that if $\xi = x - s$, then
\[ \frac{\pd}{\pd x}G = \frac{d}{d\xi}G \frac{\pd \xi}{\pd x} = G' \]
and
\[ \frac{\pd}{\pd s}G = \frac{d}{d\xi}G \frac{\pd \xi}{\pd s} = -G' \]
were used in the above derivation. So, under the assumptions made above, 
the derivative and filtering operations are interchangeable in the case 
of uniform filter width which occurs in homogeneous turbulence.

For the non--homogeneous wall bounded flow, since the scale of the large and
small eddies will depend on position relative to the wall, it is clear that
the filter width must be a function of location. To emphasize this, we rewrite
the definition of filtering, Eqn.~(\ref{defn-filter}), as 
\begin{equation}
\label{defn-filter-1}
  \ol{f}({\bf x}) = \frac{1}{\del(\bx)}
     \int_\Gamma G\left(\frac{\bx-\bs}{\del(\bx)}\right) f(\bs) \,d\bs.
\end{equation}
It can be shown by simply carrying out the differentiation, that the derivative
and filtering operations are not, in general, interchangeable. For this
reason, there is ongoing research to find filters that reduce the error
in interchanging the derivative and filtering operations 
(Germano \cite{germano3}; Ghosal and Moin \cite{ghosal1}; van der Ven
\cite{ven}; Najjar and Tafti \cite{najjar}; 
Vasilyev and Lund \cite{vasilyev} \cite{vasilyev1}; and Marsden, Vasilyev and
Moin \cite{marsden}.) 
However, as of the time of this report, 
there is no compelling evidence to indicate the best filter choice. As
mentioned before, in using the finite difference methods, researchers often
assumed an implicit filter. The issue of the interchange of filter and 
differentiation was then bypassed by assuming that just as when using the
Reynolds Average method, one assumes that all turbulent behavior is captured
by the model for the Reynolds stress, all of the effects of filtering
were assumed to be captured by the sub--grid scale model. More recently, there
have been some studies where an explicit discrete filter is used instead of
the implicit filter (Najjar et~al.~\cite{najjar};
Vasilyev et~al.~\cite{vasilyev} \cite{vasilyev1}; and
Marsden et~al.~\cite{marsden}.) However, as Vasilyev pointed out, there
is more work to be done before any conclusions can be drawn about the
effectiveness of the discrete filter.

In the case of the Top--Hat filter, in 1-dimension, the filter of $f$ is
\[ \ol{f} = \frac{1}{\del(x)} \int^{x+\frac{\del(x)}{2}}_{x-\frac{\del(x)}{2}}
            f(s) \,ds \]
and
\[ \ol{\frac{\pd f}{\pd x}} = 
       \frac{1}{\del(x)} \int^{x+\frac{\del(x)}{2}}_{x-\frac{\del(x)}{2}}
        f'(s) \,ds
    = \frac{f(x+\frac{\del(x)}{2}) - f(x-\frac{\del(x)}{2})}{\del(x)}. \]
Using the above, we calculate the derivative of the filter of $f$ as follows:
\begin{eqnarray*}
  \frac{\pd \ol{f}}{\pd x} 
&=& \frac{\pd}{\pd x} \left[
    \frac{1}{\del(x)} \int^{x+\frac{\del(x)}{2}}_{x-\frac{\del(x)}{2}}
            f(s) \,ds \right] \\
&=& -\frac{\del'(x)}{\del^2(x)} 
    \int^{x+\frac{\del(x)}{2}}_{x-\frac{\del(x)}{2}} f(s) \,ds
  + \frac{1}{\del(x)} 
    \left[ 
          \left(1+\frac{\del'(x)}{2} \right) f\left(x+\frac{\del(x)}{2}\right)
    \right] \\
&&
  - \frac{1}{\del(x)}
    \left[ 
          \left(1-\frac{\del'(x)}{2} \right) f\left(x-\frac{\del(x)}{2}\right)
    \right] \\
&=& \ol{\frac{\pd f}{\pd x}} - \frac{\del'(x)}{\del(x)} \ol{f}
  + \frac{\del'(x)}{2\del(x)}
    \left[ 
      f\left(x+\frac{\del(x)}{2}\right) + f\left(x-\frac{\del(x)}{2}\right)
    \right]. \\
\end{eqnarray*}
By the trapezoidal rule, 
\[
 \ol{f}(x) 
  = \frac{1}{\del(x)} \int^{x+\frac{\del(x)}{2}}_{x-\frac{\del(x)}{2}}f(s) \,ds
  = \frac{1}{2}
    \left[ 
      f\left(x+\frac{\del(x)}{2}\right) + f\left(x-\frac{\del(x)}{2}\right)
    \right]
  + O(\del^2 f'').
\]
So, using this and the previous equation, we get
\[ 
  \frac{\pd \ol{f}}{\pd x} - \ol{\frac{\pd f}{\pd x}} =  
  - \frac{\del'(x)}{\del(x)} \ol{f}
  + \frac{\del'(x)}{\del(x)} [\ol{f} - O(\del^2 f'')]
  = O(\del' \del f'').
\]
Since we choose our filter width, $\del$, to be the length of our grid
element and the grids that we use are graded near the wall, $\del' \not= 0$
and so the above shows that the magnatude of the error in the interchange 
of the derivative and filtering operator is dependent on the derivative of
the filter width, $\del'$. To obtain a better prespective on the implications
of the interchage error, we first rewrite it as
\[  
   \frac{\pd \ol{f}}{\pd x} - \ol{\frac{\pd f}{\pd x}} = 
     O(\frac{\del'}{\del}\del^2 f'').
\]
It is clear that if $\frac{\del'}{\del} = O(1)$, then the interchage error
is second order. Now, assuming that
\[ \frac{\del'}{\del} = M \]
for some constant, $M$, we get
\[ \del = Ce^{Mx} \]
where $C$ is a constant. Hence, if the filter width behaves like the 
exponential function, then the error in the interchange of filtering and
differentiation is second order, using the Top--Hat filter.

\subsection{Filtering and FEM}
\label{fltr-fem-sec}

One of the advantages and/or characteristics of the Finite Element Method (FEM) 
over the Finite Difference Method is that instead of producing a solution over 
a discrete set of points, it produces a continuous approximation of the actual
solution. In fact the general shape of the approximation function is known and
chosen by the user via the shape functions before the actual computations
begin. This knowledge of the general shape of the function being filtered 
naturally leads to the question of how this information can be used or effects
the analysis of the error in the interchange of filtering and differentiation.

We begin by modifying the definition of filtering, for the one dimensional case,
in the following way.
\[
   \ol{u}(x) = \frac{1}{\del(x)}
              \int_{x-\frac{\del(x)}{2}}^{x+\frac{\del(x)}{2}} 
                    G\left(\frac{x-s}{\del(x)}\right) u_e(s) \,ds
\]
where $u_e(x)$ is the finite element representation of $u(x)$.
Now using the transformation $y = \frac{x-s}{\del(x)}$, the above becomes
\[
   \ol{u}(x) = \int^{\frac{1}{2}}_{-\frac{1}{2}}
               G(y)u_e(x-y\del(x)) \,dy
\]
Differentiating this, we get
\begin{equation}
\label{fem-fltr-err1}
  \frac{d\ol{u}}{dx}(x) = \int^{\frac{1}{2}}_{-\frac{1}{2}}
               G(y)u'_e(x-y\del(x))(1-y\del'(x)) \,dy.
\end{equation}
Note that in general $u_e(x)$ is only required to be continuous rather than
continuously differentiable. It is however continuously differentiable in each
element. So we use the above notation with the understanding that
\begin{eqnarray*}
   \int^{\frac{1}{2}}_{-\frac{1}{2}}
               G(y)u'_e(x-y\del(x))(1-y\del'(x)) \,dy
   &=& \int_{-\frac{1}{2}}^{x_1} G(y)u'_{e_1}(x-y\del(x))(1-y\del'(x)) \,dy
   \\
   &+& \int_{x_1}^{x_2} G(y)u'_{e_2}(x-y\del(x))(1-y\del'(x)) \,dy
   \\
   &+& \cdots
   \\
   &+& \int_{x_n}^{\frac{1}{2}} G(y)u'_{e_n}(x-y\del(x))(1-y\del'(x)) \,dy
\end{eqnarray*}
where $(x_i, x_{i+1})$ represents the transformed elements contained in the
support for the filter.
Note that
\[
   \ol{\frac{du}{dx}}(x) = \int^{\frac{1}{2}}_{-\frac{1}{2}}
               G(y)u'_e(x-y\del(x)) \,dy
\]
and so using this in Eqn.~(\ref{fem-fltr-err1}) we get our error term
\begin{equation}
\label{fem-deriv-err}
   \left|\frac{d\ol{u}}{dx}(x)-\ol{\frac{du}{dx}}(x)\right| =
     \left|\del'(x)\int^{\frac{1}{2}}_{-\frac{1}{2}}
               yG(y)u'_e(x-y\del(x)) \,dy.\right|
\end{equation}

If we assume a linear shape function then
\[ u_{e_i}(x) = \alpha_i + \beta_i(x) \]
for $x \in \Omega_i$ where $\Omega_i$ is the i-th element. We define the filter
width $\del(x)$ to be twice the width of the smaller of the two intervals 
containing the nodal point $x$. This allows us to seperate the error term,
Eqn.~(\ref{fem-deriv-err}), in terms of each element as follows:
\[
   \left|\frac{d\ol{u}}{dx}(x)-\ol{\frac{du}{dx}}(x)\right| =
     \beta_1 \del'(x) \int^{\frac{1}{2}}_0 yG(y) \,dy
     + \beta_2 \del'(x) \int^0_{-\frac{1}{2}} yG(y) \,dy.
\]
The above would seem to indicate that the error in the interchange is of order
$\del'(x)$. This would be the case if the filter width is defined such that 
a change in the filter width does not imply a change in the element width.
However, the above definition of the filter width not only involves
$x$ but also the element width, which we denote $\del_{e_i}(x)$. Now, we note 
that
\[
   \beta_1 = \frac{u(x) - u(x-\del_{e_1}(x))}{\del_{e_1}(x)} 
           = u'(x) + O(\del_{e_1}(x))
\]
and 
\[
   \beta_2 = u'(x) + O(\del_{e_2}(x))
\]
where if we let $x= x_n$, then $\del_{e_1}(x) = x_n - x_{n-1}$ and 
$\del_{e_2}(x) = x_{n+1} - x_n$ for gradings of $x_n$ that increase with $n$.
Since $G$ is defined as an even function we have 
\[
   \int^{\frac{1}{2}}_0 yG(y) \,dy = -\int^0_{-\frac{1}{2}} yG(y) \,dy.
\]
So
\begin{eqnarray*}
   \left|\frac{d\ol{u}}{dx}(x)-\ol{\frac{du}{dx}}(x)\right| &=&
       \left|(\beta_2-\beta_1)\del'(x)\int^{\frac{1}{2}}_0 yG(y) \,dy\right|
       \\ 
       &=& \Bigg|(O(\del_{e_2}(x)) - O(\del_{e_1}(x)))
               \del'(x)\int^{\frac{1}{2}}_0 yG(y) \,dy\Bigg|
       \\
       &=& O(\del_e(x) \del'(x))
\end{eqnarray*}
where for practical problems the order of the local element widths are the same,
i.~e. $O(\del_{e_1}(x)) = O(\del_{e_2}(x)) = O(\del_e(x))$. Furthermore, if
we assume a smooth exponential grading such that $O(\del'(x))=O(\del(x))$, then
the interchange error is second order,
\[
   \left|\frac{d\ol{u}}{dx}(x)-\ol{\frac{du}{dx}}(x)\right|
     = O(\del^2_e(x)).
\]

For a quadratic shape function, we have
\[
   u_{e_i}(x) = \alpha_i + \beta_i x + \gamma_i x^2
\]
for $x \in \Omega_i$. Again, using the above filter definition, the 
interchange error is
\begin{eqnarray*}
   \left|\frac{d\ol{u}}{dx}(x)-\ol{\frac{du}{dx}}(x)\right| &=&
     (\beta_1+2\gamma_1 x) \del'(x) \int^{\frac{1}{2}}_0 yG(y) \,dy \\
     &-& 2\gamma_1\del(x)\del'(x) \int^{\frac{1}{2}}_0 y^2G(y) \,dy \\
     &+& (\beta_2+2\gamma_2 x) \del'(x) \int^0_{-\frac{1}{2}} yG(y) \,dy \\
     &-& 2\gamma_2\del(x)\del'(x) \int^0_{-\frac{1}{2}} y^2G(y) \,dy
\end{eqnarray*}
Note that in FEM, $u_{e_i}(x) = \alpha_i + \beta_i x + \gamma_i x^2$ is exact 
at each nodal point, i.~e. $u_{e_i}(x_n) = u(x_n)$. 
So, using the above definition of $\del_{e_i}(x)$, we 
have the following first order approximation of the derivative of $u$:
\begin{eqnarray*}
  u'(x) &=& \frac{u(x) - u(x-\del_{e_1}(x))}{\del_{e_1}(x)} + O(\del_{e_1}(x))
    \\
    &=& \frac{u_{e_1}(x) - u_{e_1}(x-\del_{e_1}(x))}{\del_{e_1}(x)} 
              + O(\del_{e_1}(x))
    \\
    &=&
    \frac{\alpha_1 + \beta_1 x + \gamma_1 x^2
      - [\alpha_1 + \beta_1(x-\del_{e_1}(x)) + \gamma_1(x-\del_{e_1}(x))^2]}
         {\del_{e_1}(x)} + O(\del_{e_1}(x))
    \\
    &=& \beta_1 + 2\gamma_1 x - \gamma_1\del_{e_1}(x) + O(\del_{e_1}(x))
\end{eqnarray*}
So
\[
   \beta_1 + 2\gamma_1 x = u'(x) + O(\del_{e_1}(x))
\]
and using a similiar argument, we also have
\[
   \beta_2 + 2\gamma_2 x = u'(x) + O(\del_{e_2}(x)).
\]
Hence, using the symmetry of $G$ as in the linear case, the interchange 
error is again
\[
   \left|\frac{d\ol{u}}{dx}(x)-\ol{\frac{du}{dx}}(x)\right| =
   O(\del_e(x)\del'(x)).
\]

We now show that for any polynomial shape function, the interchange error is
of order $O(\del_e(x)\del'(x))$. Consider a general n-th order polynomial 
shape function
\[
   u_{e_j}(x) = \sum_{i=0}^{n} \alpha_{i,j} x^i.
\]
Note that 
\[
   u'_{e_j}(x) = \sum_{i=1}^n i\alpha_{i,j} x^{i-1}.
\]
So the interchange error, Eqn.~(\ref{fem-deriv-err}), can be written as
\begin{eqnarray*}
   \left|\frac{d\ol{u}}{dx}(x)-\ol{\frac{du}{dx}}(x)\right| &=&
       \Bigg|\del'(x)\int^0_{-\frac{1}{2}}
               yG(y)\sum_{i=1}^n i\alpha_{i,1}(x-y\del(x))^{i-1} \,dy
	\\
	 &+& \del'(x)\int^{\frac{1}{2}}_0
               yG(y)\sum_{i=1}^n i\alpha_{i,2}(x-y\del(x))^{i-1} \,dy \Bigg|.
\end{eqnarray*}
We can ignore any term that has $\del(x)$ since they will be of oder
$O(\del(x)\del'(x))$ or higher. So, using the Binomial Theorem and the fact
that $G$ is symmetric, the 
interchange error above can be rewritten as
\begin{eqnarray}
\label{fem-general-fltr-err}
   \left|\frac{d\ol{u}}{dx}(x)-\ol{\frac{du}{dx}}(x)\right| &=&
       \Bigg|\del'(x)
       \Big[(\sum_{i=1}^n i\alpha_{i,2} x^{i-1})
               -(\sum_{i=1}^n i\alpha_{i,1} x^{i-1})\Big]
          \int^{\frac{1}{2}}_0 yG(y) \,dy
\\ \nonumber
	&+& O(\del(x)\del'(x))\Bigg|
\end{eqnarray}
Now, as in the quadratic case, we note that the interpolating polynomial
function is exact at each nodal point, i.~e. $u_{e_j}(x_n)=u(x_n)$. So
\begin{eqnarray*}
   u'(x) + O(\del_{e_1}(x)) &=& \frac{u(x) - u(x-\del_{e_1}(x))}{\del_{e_1}(x)}
   \\
   &=& \frac{u_{e_1}(x) - u_{e_1}(x-\del_{e_1}(x))}{\del_{e_1}(x)} 
   \\
   &=& \frac{\sum_{i=0}^n \alpha_{i,1}x^i 
           - \sum_{i=0}^n \alpha_{i,1}(x-\del_{e_1})^i(x)}{\del_{e_1}(x)}
\end{eqnarray*}
Again, using the Binomial Theorem, we rewrite the above as
\[
   u'(x) + O(\del_{e_1}(x)) 
     = \frac{\sum_{i=0}^n \alpha_{i,1}x^i 
           - \sum_{i=0}^n \alpha_{i,1}
             \sum_{j=0}^i{i \choose j}x^j(-\del_{e_1}(x))^{i-j}}
            {\del_{e_1}(x)} 
\]
Note that the terms that do not have any factor of $\del_{e_1}(x)$ cancel 
out which leaves the terms with a factor of just $\del_{e_1}(x)$ as the lowest
order terms. So, collecting these terms, the above can be rewritten as
\[
   u'(x) + O(\del_{e_1}(x)) = \sum_{i=1}^n i\alpha_{i,1} x^{i-1}.
\]
Using same argument as above only with a first order forward differencing 
scheme instead of the backward differencing scheme, we get
\[ 
   u'(x) + O(\del_{e_2}(x)) = \sum_{i=1}^n i\alpha_{i,2} x^{i-1}.
\]
Using the above in Eqn.~(\ref{fem-general-fltr-err}), we get
\[
   \left|\frac{d\ol{u}}{dx}(x)-\ol{\frac{du}{dx}}(x)\right|
    = O(\del_e(x)\del'(x))
\]
where again we assume that the filter and element widths are
of the same order, i.~e. $O(\del(x)) = O(\del_{e_i}(x)) = O(\del_e(x))$.
So if we assume that the grading of the elements is smooth enough such that
$O(\del'(x))=O(\del(x))$ then the error in the interchange of differentiation
and filtering is second order for any polynomial interpolation of degree one
or greater. 

Appendix \ref{fltr-err-comp} presents the results of a computational
experiment that verifies the above results for the case of linear and 
quadratic shape functions. Note that since the FEM is applied to solve the
filtered equations, the above filtering technically only applies to filtering
a filtered quantity. This is true for any explicit filter of $u$ since the
only values of $u$ that are available are those that are obtained from the
numercial method. But the numerical method is used to solve the LES or filtered
equations and not the original, as discussed in the following section. Hence,
to be completely consistent, an explicit filter can only be applied in a 
double filtering situation such as when using the Dynamic Subgrid Scale model
discussed below. It would seem that to be consistent the initial filtering
needs to be introduced in one of two ways. One, it can be introduced by a term
representing the derivative/filter interchange error or this error can be
ignored by using a filter whose interchange error is of the same order as the
numerical scheme. Note that this error must be derived using the original
vector field $u$ and should avoid any reference to the computed $u$.
Also, this filter would not be applied directly to the
numerical solution but should instead be applied to the experimental or 
DNS data for comparison with the LES data. The second and ideal place to 
introduce the filter effect would be the SGS model since the filter determines
the scales that are in the system.

\subsection{Governing Equations}

The set of equations that govern the motion of (Newtonian) fluids are known
collectively as the Navier--Stokes equations. In the case of incompressible
fluids with uniform density, the non--dimensionalized equations are
\[ \left\{
  \begin{array}{l}
    \frac{\pd u_\alpha}{\pd x_\alpha} = 0 \\
    \frac{\pd u_\alpha}{\pd t} 
    + \frac{\pd}{\pd x_\beta}(u_\alpha u_\beta)
    =
    - \frac{\pd P}{\pd x_\alpha} 
    + \frac{1}{Re} \frac{\pd^2}{\pd x^2_\beta} u_\alpha ,
  \end{array}
  \right. 
\]
where we use the standard tensor notation of repeated indices indicating
summation, $P = \frac{P}{\rho}$, and $Re = \frac{UL}{\nu}$ with $U$ begin
the characteristic velocity, $L$ begin the characteristic length, and
$\nu$ being the kinematic viscosity. 
Now, assuming that the error in the interchange of filtering and differentiation
can be ignored (i.~e.~$\frac{\ol{\pd f}}{\pd x} = \frac{\pd \ol{f}}{\pd x}$),
the filtered Navier--Stokes equations are
\begin{equation}
  \label{les-NS}
  \left\{
  \begin{array}{l}
    \frac{\pd \ol{u}_\alpha}{\pd x_\alpha} = 0 \\
    \frac{\pd \ol{u}_\alpha}{\pd t} 
    + \frac{\pd}{\pd x_\beta}(\ol{u_\alpha u_\beta})
    =
    -\frac{\pd \ol{P}}{\pd x_\alpha} 
    + \frac{1}{Re} \frac{\pd^2}{\pd x^2_\beta} \ol{u}_\alpha.
  \end{array}
  \right.
\end{equation}
By filtering the equations above, we have introduced a new unknown term,
$\ol{u_\alpha u_\beta}$, but no new equation. This means that the system of 
equations
is now underdetermined. This is known as the closure problem, which is also
encountered by researchers using RAS. In fact, 
the Reynolds Averaged Navier--Stokes equations look very much like the 
filtered Navier--Stokes equation above (Eqn.~(\ref{les-NS}) with the 
$\ol{\cdot}$ replaced by $<\cdot>$, which represents the Reynolds averaging
operator.) The difference is in how the 
$\ol{\cdot}$ (and/or $<\cdot>$) is defined. In LES, the decomposition 
of variables
is in terms of the large scales and small, i.~e.~$u = \ol{u} + u'$ where
$\ol{u}$ is the large scale component of $u$ and $u'$ is the small or
sub--scale component. From the definition of filtering, 
Eqn.~(\ref{defn-filter}), it is clear that in general 
$\ol{\ol{u}} \not= \ol{u}$, unless of course $u$ is constant or the filter
used is the Sharp--Cutoff. This means that,
in general, $\ol{u'} \not= 0$. In the Reynolds Averaging method, the 
decomposition is in terms of the statistical mean and fluctuations with 0 mean. So if we write $u = <u> + u'$ then $<<u>> = <u>$ and $<u'> = 0$. This 
difference, in the interpretation of the decomposition of  
variables, leads to a very different interpretation, both physically and 
mathematically, of Eqn.~(\ref{les-NS}). In either case, note that the
interaction between large and small scales ($\ol{u}$ and $u'$) or mean
and fluctuations ($<u>$ and $u'$) must occur in the $\ol{u_\alpha u_\beta}$ 
term of Eqn.~(\ref{les-NS}). As this work is based on LES, our focus will be
on the closure problem as it relates to LES. 

To deal with the closure problem, the above filtered equations,
Eqn.~(\ref{les-NS}), are rewritten as 
\begin{equation}
  \label{les-NS-1}
  \left\{
  \begin{array}{l}
    \frac{\pd \ol{u}_\alpha}{\pd x_\alpha} = 0 \\
    \frac{\pd \ol{u}_\alpha}{\pd t} 
    + \frac{\pd}{\pd x_\beta}(\ol{u}_\alpha \ol{u}_\beta)
    =
    - \frac{\pd \ol{P}}{\pd x_\alpha} 
    - \frac{\pd}{\pd x_\beta}\tau_{\alpha \beta}
    + \frac{1}{Re} \frac{\pd^2}{\pd x^2_\beta} \ol{u}_\alpha
  \end{array}
  \right.
\end{equation}
where 
\[ \tau_{\alpha \beta} = \ol{u_\alpha u_\beta} - \ol{u}_\alpha \ol{u}_\beta \]
is called the subgrid scale (SGS) Reynolds stress. The SGS Reynolds stress 
in LES is similiar to the Reynolds stress in RAS in that in both cases the 
the unresolved scales, $u'$ (small eddies in LES and turbulent fluctuations 
in RAS,) are viewed as producing stresses in the resolved scales (large eddies
in LES and the mean velocity in RAS.) The difference is that in RAS, $u'$
represents all the turbulent motions, while in LES, $u'$ represents only
the small eddies or subgrid scales. This means that the energy in the 
unresolved scales of LES (subgrid scales) is a much smaller portion of the 
total flow compared to the energy in the unresolved scales of RAS (the 
turbulent fluctuations.) Hence, it is believed that the accuracy in the
modelling of the SGS stress is not as crucial as the modelling of the Reynolds
stress term in RAS. 

It should be noted here that our treatment of the closure problem, 
Eqn.~(\ref{les-NS-1}), is what Mason \cite{mason} refers to as the 
Lilly-Deardorff approach. This approach allows us to define the filtering
operator implicitly through the modelling of the SGS stress term.
Another approach, referred to as the Leonard approach, is to substitute 
\[
   u_\alpha = \ol{u}_\alpha + u'_\alpha 
\]
into $\ol{u_\alpha u_\beta}$ to get
\[ 
  \ol{u_\alpha u_\beta} = \ol{\ol{u}_\alpha\ol{u}_\beta} 
             + \ol{u'_\alpha \ol{u}_\beta} + \ol{\ol{u}_\alpha u'_\beta}
             + \ol{u'_\alpha u'_\beta}.
\]
Note that the first term on the right hand side of the above equation can be
explicitly calculated if the filter used is itself explicitly defined. Using
the above, the SGS stress term in Eqn.~(\ref{les-NS-1}) becomes
\begin{equation}
  \label{l-stree}
  \tau_{\alpha \beta} 
      = (\ol{\ol{u}_\alpha\ol{u}_\beta} - \ol{u}_\alpha \ol{u}_\beta)
      + (\ol{u'_\alpha \ol{u}_\beta} + \ol{\ol{u}_\alpha u'_\beta})
      + \ol{u'_\alpha u'_\beta}.
\end{equation}
The advantage of the Leonard approach is that it breaks the SGS stress 
into three terms that have physical interpretations in the flow.
The first term on the right hand side of Eqn.~(\ref{l-stree}) represents
the interaction of large eddies that produce small eddy effects and is
called the Leonard stress. The second term, called the cross term, represents
the interaction between the large and small eddies and the third term
represents the interaction of small eddies to produce large eddy effects. It
is this term that produces the transfer of energy from the small to the
large eddies and soe is known as the backscatter term. Note that the second
term may also produce backscatter. The disadvantage of the Leonard apprach
is that the filter must be explicitly known for the calculation of the Leonard
stress term. This may not be a trivial task, especially in grid base methods
where the filter is often implicitly defined. As pointed out by 
Mason~\cite{mason} and Ferziger~\cite{ferziger}, both the Lilly-Deardorff
and Leonard approaches have been used successfully by researchers. The 
approach taken for our work is that of Lilly-Deardorff.

\subsection{Subgrid Scale Models}

As discussed in the previous section, in filtering the Navier--Stokes 
equations a new term is introduced which leads to a closure problem. To 
deal with this problem, the filtered Navier--Stokes equations are rewritten
into the form given in Eqn.~(\ref{les-NS-1}) and the new term,
$\tau_{\alpha \beta}$, is then modelled. In this section, we will review two
of the most commonly used SGS models, namely the Smagorinsky and the
Dynamic Subgrid Scale models, and introduce our variant of the Smagorinsky
model, which we will refer to as the Modified Smagorinsky model. 

\subsubsection{Smagorinksy Model}

The Smagorinsky \cite{smag} subgrid scale model,
\begin{equation}
  \label{smag-stress}
  \tau_{\alpha \beta} - \frac{1}{3}\tau_{\gamma \gamma} \delta_{\alpha \beta}
    = -\nu_T
       \left( 
              \frac{\pd \ol{\ua}}{\pd \xb} + \frac{\pd \ol{\ub}}{\pd \xa}
       \right)
    = -2 \nu_T \ol{S}_{\alpha\beta}
\end{equation}
where $\nu_T$ is called the eddy viscosity (to be derived below) and 
\[
  \ol{S}_{\alpha \beta} = \frac{1}{2}
                    \left( 
                           \frac{\pd \ol{\ua}}{\pd \xb} 
                         + \frac{\pd \ol{\ub}}{\pd \xa}
                    \right),
\]
is known as an 
eddy viscosity type model because it assumes that the small eddies remove
energy from the flow through a dissipative process. This precludes the 
possiblity of backscatter (the reintroduction of energy from the small to
large eddies) which is considered one of the model's weak points. The 
following derivation of the model is that of Ferziger \cite{ferziger}.

According to turbulence theory, energy is introduced at the largest scales
and is successively transferred to the smaller scales until viscous damping
becomes the dominant effect. In the region where the viscous effect becomes
dominant, the turbulent energy of the flow is damped out by the transfer of 
the kinetic energy to internal energy. Between these two scales is a region
known as the inertial subrange, where there is neither significant production
nor dissipation of energy. In this region, only the inviscid mechanisms are
active and so assuming that the transfer of energy is always from the large
to small scales, the term responsible for the transfer of this energy is then
the advection term of the Navier--Stokes equations. This allows us to estimate
the rate of energy transfer to the small scales as the magnitude of the 
contribution of the advection term to the kinetic energy equation, which is
\[ 
  \frac{1}{2} \frac{\pd}{\pd \xb}(\ua \ua \ub).
\]
Since energy is introduced at the largest scales, we estimate the dissipation
rate, $\epsilon$, as 
\begin{equation}
  \label{diss-large}
  \epsilon \approx U_L^3/L
\end{equation}
where $U_L$ is the velocity scale of the large eddies and $L$ is the integral
length scale of the turbulent flow. Furthermore, if we assume that the largest
subgrid scales are much larger than the viscous scales, then
\begin{equation}
  \label{diss-small}
  \epsilon \approx U_{SGS}^3/\del
\end{equation}
where $U_{SGS}$ is the velocity scale of the small eddies and $\del$ is the
length of the largest subgrid scale eddies, which is also the length scale
associated with the filter. 

Now, as stated before, turbulence theory states that energy is transferred
from the large to small scales. We assume that this process is dissipative
in nature with respect to the large eddies (i.~e.~once the energy is lost
by the large eddies to the small, it cannot be recovered.) So, the SGS
model represents this energy transfer as effective viscous dissipation. 
Since it is the smallest resolved scales (of size $\del$) that are most
influenced by the SGS model, we estimate the effective dissipation as
\begin{equation}
  \label{diss-model}
  \epsilon \approx \nu_T U_{SGS}^2/\del^2.
\end{equation}
Using (\ref{diss-small}) and (\ref{diss-model}), we get
\[ 
  \nu_T \propto \del U_{SGS}
\]
and using the above together with Eqns.~(\ref{diss-large}) and 
(\ref{diss-small})
\[
  \nu_T \approx \del^{\frac{4}{3}} L^{-\frac{1}{3}} U_L.
\]
Using
\[
  U_L \approx L(2\ol{S}_{\alpha \beta} \ol{S}_{\alpha \beta})^{\frac{1}{2}}
      = L |\ol{S}|
\]
we get
\[
  \nu_T = C_S^2 \del^{\frac{4}{3}} L^{\frac{2}{3}} |\ol{S}|
\]
where $C_S$ is the Smagorinsky parameter which is introduced to produce
the equality. Because of the difficulty in calculating the integral length
scale, $L$, most researchers make the simplifying substitution
\[ 
  \del^{\frac{4}{3}} L^{\frac{2}{3}} \to \del^2
\]
which leads to the more common form of the eddy viscosity in the
Smagorinsky model,
\begin{equation}
  \label{eddy-vis}
   \nu_T = (C_S \del)^2 |\ol{S}|.
\end{equation}
In our calculations we use as our length scale
\[ 
  \del = (\del_1 \del_2 \del_3)^{\frac{1}{3}}
\]
where $\del_\alpha$ is the grid element length in the appropriate direction
and set 
\[
  C_S = 0.065.
\]

Although the Smagorinsky model is the most widely recognized and used SGS
model, it is not without its problems. One significant problem is that
in non--homogeneous wall bounded flows, near the surface of the walls, the
model does not damp out the eddy viscosity enough to allow kinematic viscosity
effects to become dominant. In fact, as the wall and hence the viscous sublayer
is approached, our approximation, (\ref{diss-small}), in the derivation of the
model itself becomes questionable since the grid points, hence $\del$, will
be in or near the viscous sublayer region. 

Another problem as mentioned before, is that since the Smagorinsky model 
assumes that the transfer of energy from large to small scales is a 
diffusive process, Eqn.~(\ref{diss-model}), no mechanism exists to 
transfer energy from the small to
the large scales, i.~e.~there is no backscatter. To show this, we follow
the work of Piomelli $et~al.$~\cite{piomelli2} and observe that in the
resolved energy transport equation,
\[
  \frac{\pd \ol{q}^2}{\pd t} + \frac{\pd(\ol{q}^2\ol{u}_\beta)}{\pd \xb}
  = \frac{\pd}{\pd \xb} \left( -2\ol{p} \, \ol{u}_\beta 
  - 2\ol{u}_\alpha \tau_{\alpha \beta}
  + \frac{1}{Re} \frac{\pd \ol{q}^2}{\pd \xb} \right)
  - \frac{2}{Re} \frac{\pd \ol{u}_\alpha}{\pd \xb} 
    \frac{\pd \ol{u}_\alpha}{\pd \xb}
  + 2 \tau_{\alpha \beta} \ol{S}_{\alpha \beta}.
\]
where $\ol{q}^2 = \ol{u}_\alpha \ol{u}_\alpha$, one--half of
the last term,
\[ 
  \epsilon_{SGS} = \tau_{\alpha \beta}\ol{S}_{\alpha \beta}
\]
which is referred to as the ``subgrid--scale dissipation'',
represents the transfer of energy between the large and small scales. Note
that 
\[
  \epsilon_{SGS} < 0
\]
indicates that the resolved energy decays with respect to the interaction
between the large and small scales (forward scatter), while 
\[
  \epsilon_{SGS} > 0
\]
means that the resolved energy grows with respect to the interaction of
the large and small scales (backscatter). In the case of the Smagorinsky model,
Equations (\ref{smag-stress}) and (\ref{eddy-vis}) cleary show that 
$\epsilon_{SGS} < 0$, and so there is no backscatter when the Smagorinsky
model is used. In their study using DNS data, Piomelli 
$et~al.$~\cite{piomelli2} found that
for channel flow using a Sharp--Cutoff filter, nearly 50\% of the points
experienced backscatter. For the Gaussian filter, the percentage dropped to 
about 30\% and the Top--Hat filter fell inbetween the cutoff and the Gaussian.
These results indicate that although in the mean, the net effect of the SGS
dissipation is forward scatter, there is significant backscatter. The 
importance of backscatter in the SGS model is not yet clear. In fully 
developed channel flows for example, purely dissipative SGS stress models
have been used successfully. Since forward scatter is dominant in the mean,
this would seem reasonable. As Piomelli et al.~\cite{piomelli2} speculated,
it would seem reasonable that backscatter effects would be more important
for nonequilibrium flows.

Finally, there is the problem of setting the Smagorinsky parameter, $C_S$.
The values used range from about $0.2$ for isotropic turbulence to $0.065$
for the channel flow. Other than for isotropic turbulence, the author is
not aware of any systematic analysis to determine the parameter.

\subsubsection{Dynamic Subgrid Scale Model}

To overcome the problems of the Smagorinsky model, Germano 
$et~al.$~\cite{germano1} introduced a method of calculating the Smagorinsky
parameter interactively during run time rather than setting the parameter 
to a fixed constant at the beginning of the numercial calculations. 
Here, we introduce a variant of Germano's dynamic procedure developed by
Piomelli and Liu \cite{piomelli3}.
The procedure involves applying a second filter, which is called the 
test filter $\wdtld{\cdot}$, to the filtered momentum equation, 
Eqn.~(\ref{les-NS}). This results in
\[
  \frac{\pd \wdtld{\ol{u}}_\alpha}{\pd t} 
  + \frac{\pd}{\pd x_\beta}(\wdtld{\ol{u}_\alpha} \wdtld{\ol{u}_\beta})
  =
  -\frac{\pd \wdtld{\ol{P}}}{\pd x_\alpha} 
  -\frac{\pd {\cal T}_{\alpha \beta}}{\pd x_\beta}
  + \frac{1}{Re} \frac{\pd^2}{\pd x^2_\beta} \wdtld{\ol{u}}_\alpha  
\]
where
\[ 
  {\cal T}_{\alpha \beta} = \wdtld{\ol{u_\alpha u_\beta}}
                             - \wdtld{\ol{u}}_\alpha \wdtld{\ol{u}}_\beta. 
\]
Now we define the resolved turbulent stress as
\[ 
  {\cal L}_{\alpha \beta} = \wdtld{\ol{u}_\alpha \ol{u}_\beta}
                          - \wdtld{\ol{u}}_\alpha \wdtld{\ol{u}}_\beta
\]
which represents the contribution of the small resolved scales to the 
Reynolds stresses. Using the Smagorinsky model as our SGS model, we have
at the filter and test filter level,
\begin{equation}
  \label{smag-stress-1}
  \tau_{\alpha \beta} 
    - \frac{1}{3} \tau_{\gamma \gamma} \delta_{\alpha \beta} 
  \approx
                         -2C_S \del^2 |\ol{S}| \ol{S}_{\alpha \beta} 
                       = -2C_S B_{\alpha \beta} 
\end{equation}
\begin{equation}
  \label{smag-stress-2}
  {\cal T}_{\alpha \beta} 
      - \frac{1}{3}{\cal T}{\gamma \gamma} \delta_{\alpha \beta}
    \approx
      -2C_S \wdtld{\del}^2 |\wdtld{\ol{S}}| \wdtld{\ol{S}}_{\alpha \beta}
    = -2C_S A_{\alpha \beta}
\end{equation}
where we have used $C_S$ rather than $C_S^2$ in the eddy viscosity term,
Eqn.~(\ref{eddy-vis}), and $\del$ and $\wdtld{\del}$ represents the width of
the first filter and the test filter, respectively. Now, the SGS 
stresses at the filter and test 
filter level are related to the resolved turbulent stress through an algebraic
identity known as Germano's identity,
\[
  {\cal L}_{\alpha \beta} = {\cal T}_{\alpha \beta} 
  - \wdtld{\tau}_{\alpha \beta}.
\]
Note that if we use the models for the stresses, Eqns.~(\ref{smag-stress-1})
and (\ref{smag-stress-2}), then Germano's identity is no longer satisfied
exactly. Hence, we obtain the residual equation
\begin{equation}
  \label{residual-eqn}
  E_{\alpha \beta} = {\cal L}_{\alpha \beta} 
                     + 2C_S A_{\alpha \beta} 
                     - 2 \wdtld{C_S B_{\alpha \beta}}.
\end{equation}
Piomelli's method is to minimize the square of the residual assuming that
the Smagorinsky parameter, $C_S$, is known (denoted by $C_S^*$) in the 
third term on the right side of Eqn.~(\ref{residual-eqn}). So the residual 
equation becomes
\[  
  E_{\alpha \beta} = {\cal L}_{\alpha \beta} 
                     + 2C_S A_{\alpha \beta} 
                     - 2 \wdtld{C_S^* B_{\alpha \beta}} 
\]
and minimizing the square of the residual yields
\begin{equation}
  \label{dsgs-C}
  C_S = -\frac{1}{2} \frac{({\cal L}_{\alpha \beta} 
                           - 2 \wdtld{C_S^* B_{\alpha \beta}}) A_{\alpha \beta}}
                         {A_{\xi \eta} A_{\xi \eta}}. 
\end{equation}
The Smagorinsky parameter, $C_S^*$ is then approximated using a Taylor series
expansion of $C_S$. 

One of the advantages of the DSGS model is that there are no parameters to
determine before the numerical run. Also, the DSGS model is able to reduce the
eddy viscosity in the laminar regions of the flow. Finally, the DSGS model
can allow for backscatter in the flow. 

Although Piomelli's version of the DSGS model was formulated to answer some
of the problems with Germano's original dynamic procedure, there are still
some unresolved issues with the model. One of the problems is that in the
derivation of the model, the Smagorinsky parameter was assumed to be 
independent of the filter width (see Eqns.~(\ref{smag-stress-1}) and 
(\ref{smag-stress-2}).) However, using DNS data from isotropic turbulence,
Meneveau \cite{meneveau1} found that the instantaneous values of the 
Smagorinsky parameter varied significantly when computed at two different
filter widths. 

Another problem with the DSGS model is that, as formulated,
it allows $C_S$ and hence the eddy viscosity to be negative. Although 
this allows the model to account for backscatter, if the eddy viscosity remains
negative long enough, it can lead to numerical instabilities. Researchers
have found this to be case and some type of an artifical constraint is usually
required to prevent this instability. In our case, we use the constraint
\[
  C_S \ge C_{min}
\]
where $C_{min}= -0.01$. 

Also, during the execution of the DSGS model, it was found that if we calculate
$C_S$ at each time step, then Eqn.~(\ref{dsgs-C}) can and does become
unstable. So, to stabilize the calculations, we use time averaged values to
calculate $C_S$. The time interval used is kept small (100 time steps) to 
preserve the dynamic nature of the model in time as well as in space. 

Finally, although not so much a problem as it is a question, is the issue of the
double filtering. In the derivation of the DSGS model, it is assumed that
filtering twice is equivalent to filtering once with some filter. This
assumption is used to apply the Smagorinsky SGS model to the SGS stress term
at the test filter level, Eqn.~(\ref{smag-stress-2}). In Appendix
\ref{double-filtering}, we
derive explicitly the filters that are equivalent to filtering twice with the
Gaussian and the Top--Hat filter.

\subsubsection{Modified Smagorinsky Model}

One of the problems of the Smagorinsky model as mentioned above, is that the
eddy viscosity does not damp out as the walls are approached. To overcome
this problem we introduce a damping function in the calculation of the
Smagorinsky parameter. The damping function is constructed under the
constraint that the distance from the wall cannot be used directly by the
model and that the new model cannot significantly increase resource usage
(e.~g.~cpu time and memory.) 

To construct our damping function, we first define the instantaneous
kinetic energy as
\[
  \kappa = \frac{1}{2} \ol{u}_\alpha \ol{u}_\alpha,
\]
the instantaneous dissipation rate as
\[
  \epsilon = \nu 
   \frac{\pd \ol{u}_\alpha}{\pd x_\beta}\frac{\pd \ol{u}_\alpha}{\pd x_\beta},
\]
and the instantaneous Reynolds number as
\[
  Re_{inst} = \frac{\kappa^2}{\nu \epsilon}.
\]
Note that near the wall, if we Taylor expand $\ol{u}_\alpha$ assuming that the
wall location is at $\bx = 0$, then the only non--zero derivative of the 
velocity vector will be in the wall normal direction, $y$, assuming the
wall to be in the x-z plane. So
\[
  \ol{u}_\alpha = \ol{u}_\alpha \Bigg|_{\bx=0} 
    + y \frac{\pd \ol{u}_\alpha}{\pd y} \Bigg|_{\bx=0} + \ldots
\]
By the no slip boundary condition, $\ol{u}_\alpha({\bf 0}) = 0$, and so
\[
  O(\kappa) = y^2.
\]
Using a similar argument, we have that
\[
  O(\epsilon) = 1. 
\]
So
\[ 
  Re_{inst} \sim y^4.
\]
Hence we use the instantaneous Reynolds number, $Re_{inst}$, to detect the
wall. Since the Smagorinsky model works well away from the walls, our
modified Smagorinsky model is
\begin{equation}
  \label{mod-smag}
  \nu_T = (C_S' \del)^2 |\ol{S}|
\end{equation}
where
\[
  C_S' = C_S \left\{ \begin{array}{cl}
               \left( \frac{Re_{inst}}{Re^o} \right)^\gamma 
                         & for \quad 0 \le Re_{inst} \le Re^o
\\
                1  & for\quad Re_{inst} > Re^o
\end{array}
\right.
\]
and
\[ 
  Re^o = \mu [Re_{inst}]_{max}.
\]
Note that $\gamma$ and $\mu$ are just model parameters and that 
$[Re_{inst}]_{max}$ represents the largest value of $Re_{inst}$ in the 
compuational domain at a given instant in time.

As will be shown in the results section, our model is able to detect the
walls and damp out the eddy viscosity. However, there is nothing in the
derivation of the model that indicates that $Re_{inst} \le Re^o$ only near
the wall surfaces. In fact, as will be seen in the results for the channel
and backward facing step problems, $Re_{inst} \le Re^o$ occurs quite often
in elements far from the walls.

\newpage

\section{Channel Flow}

We use the channel flow as a benchmark to test the three subgrid scale models,
because of the extensive amount of data, both experimental and theoretical, 
that is available for this type of flow. Since fully developed turbulent 
channel flow is homogeneous in both the streamwise and spanwise directions, 
we use periodic boundary conditions in those directions.

\subsection{Computational Parameters}

\begin{figure}[h]
\begin{center}
\leavevmode\epsfxsize=120mm
\leavevmode\epsfysize=80mm
\epsffile{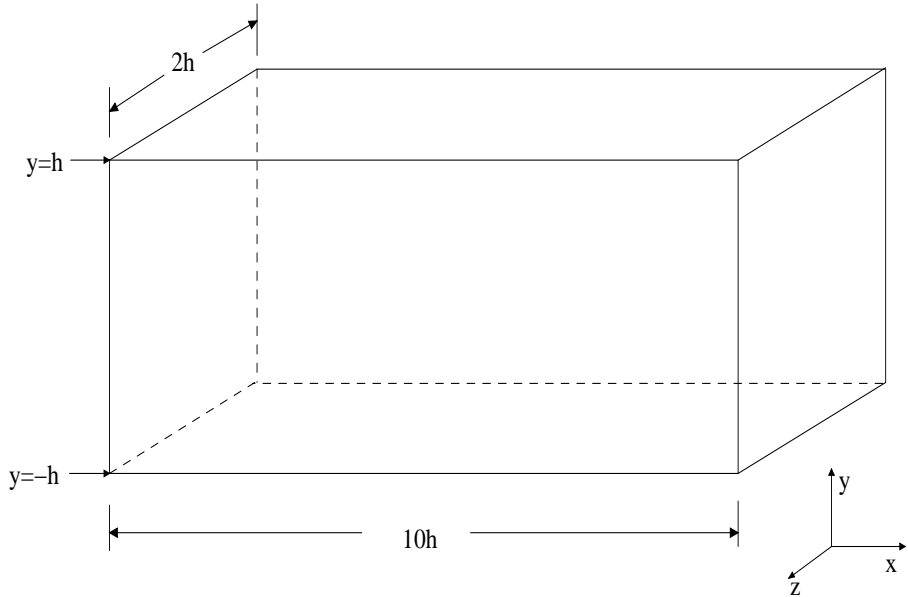}
\end{center}
\renewcommand{\baselinestretch}{1}
\caption{ \label{c-sketch}  \em
          Sketch of the channel with relevant dimensions. For our 
          experiment, we use $h=0.5$}
\end{figure}

The computational domain of the channel consists of 
$53 \times 32 \times 16$ elements and 
its dimensions are 10h in the streamwise direction (x), 2h in the spanwise 
direction (z) and 2h in the wall normal direction (y), where h $(=0.5)$ is the 
channel half height. A sketch of the channel is given in Fig.~\ref{c-sketch}.

\begin{figure}[h]
\begin{center}
\leavevmode\epsfxsize=140mm
\leavevmode\epsfysize=80mm
\epsffile{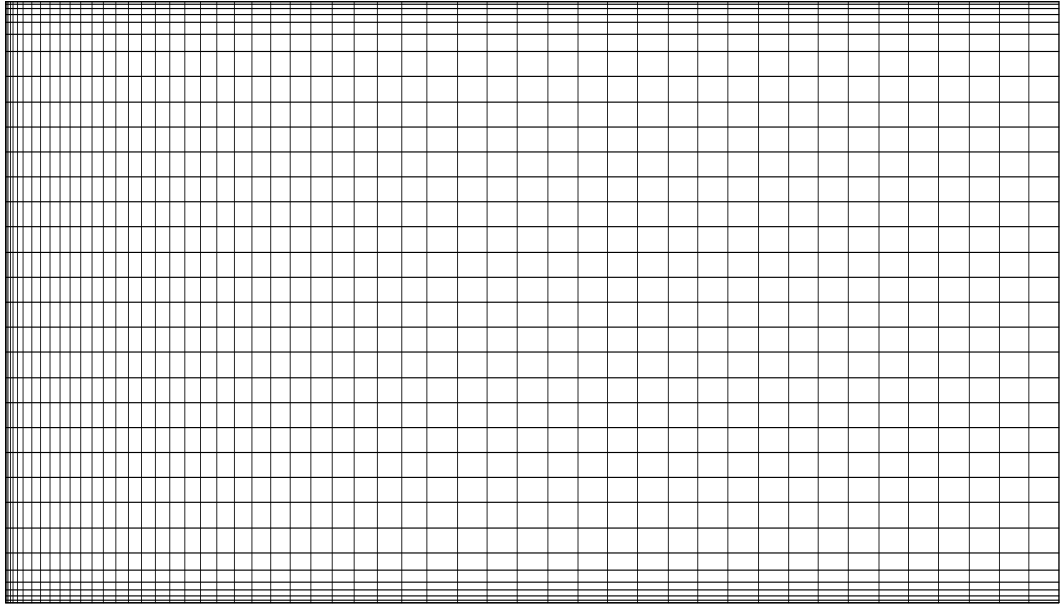}
\end{center}
\caption{ \label{xy-grid} \em
          Sketch of the xy--plane face of the graded channel grid.}
\end{figure}

For initial data, we use the following:
\[
  \left\{
  \begin{array}{l}
    u = f(y) + \epsilon u^* \\
    v = \epsilon v^* \\
    w = \epsilon w^*,
  \end{array}
  \right.
\]
where $f(y)$ is a parabolic function zero at the walls and one at the
channel center, $\epsilon = 0.1$, and $u^*,$ $v^*,$ and $w^* \in [-1, 1]$
(uniform distribution.) In our flow experiments, we set the viscosity,
$\nu = 10^{-4}.$ So, using mean initial streamwise center line velocity, 
$U_0 = 1$ as our velocity scale, and the channel height, $2h = 1,$ as our
integral length scale, we estimate the Kolmogorov length scale as
\[
  \eta \approx \left(\frac{2hU_0}{\nu} \right)^{-\frac{3}{4}} 
         = (Re_0)^{-\frac{3}{4}} = 0.001.
\]
Hence, near the walls, we grade the grid so that the smallest element will
have a height of $O(10^{-3}).$ A sketch of the graded grid in the xy--plane
face of the channel is shown in Fig.~\ref{xy-grid}. Note that the grading
in the streamwise direction was done to reduce the number of grid points 
needed in the streamwise direction. It is not a result of any physics of the
flow, as in the wall normal direction. A uniform grid is used in the spanwise
direction.

\subsection{Constant Flow Rate}

In assuming periodic boundary conditions in the streamwise direction, the
computational grid is in essence following the flow down stream. In order
to maintain a time dependent flow, an external force must be applied to the
flow. In the case of the channel, a negative streamwise pressure gradient
is used to sustain the flow. To calculate this pressure gradient, we follow
the procedure described by Tran and Morchoisne \cite{tran} in which the
pressure gradient is found in such a way that the flow rate remains 
constant. Let
\[ 
  F_\alpha = (F_x(t), 0, 0) 
\]
be the streamwise pressure gradient which may vary in time but is constant
in space. Using the centerline velocity, $U_0$, and channel half height, $h$,
as our characteristic velocity and length, the non-dimensionalized 
momentum equations are
\begin{equation}
  \label{non-dim-moment}
  \frac{\pd U_\alpha}{\pd t} + \frac{\pd}{\pd \xb}(U_\alpha U_\beta)
  = -\frac{\pd P}{\pd \xa} 
  + \frac{1}{Re_0} \frac{\pd^2 U_\alpha}{\pd x_\beta x_\beta} + F_\alpha,
\end{equation}
where $Re_0 = \frac{U_0 h}{\nu}.$

Since the flow is assumed to be homogenous in the stream and spanwise 
directions, we define our averaging operator to be
\begin{equation}
  \label{defn-avg}
  <f>_{x,z}(y) = \frac{1}{L_x L_z}
                \int_0^{L_x} \int_0^{L_z} f(x,y,z) \,dx dz,
\end{equation}
where $L_x$ and $L_z$ are the streamwise and spanwise lengths of the channel,
and the mean flow rate as
\begin{equation}
  \label{defn-Q}
  Q = \frac{L_z}{h} \int^1_{-1} <U>_{x,z} \,dy.
\end{equation}
So, integrating the streamwise momentum equation, Eqn.~(\ref{non-dim-moment}),
we get 
\begin{equation}
  \label{non-dim-moment-1}
  \frac{L_z}{h}\frac{1}{L_xL_z} \int^1_{-1} \int_0^{L_z} \int_0^{L_x}
  \left\{
    \underbrace{\frac{\pd U}{\pd t}}_A
  + \underbrace{\frac{\pd}{\pd \xb}(U U_\beta)}_B
  = \underbrace{-\frac{\pd P}{\pd x}}_C 
  + \underbrace{\frac{1}{Re_0} \frac{\pd^2 U}{\pd \xb \pd \xb}}_D 
  + \underbrace{F_x}_E
  \right\} \,d{\bf x}.
\end{equation}
We simplify each term in the above equation as follows:
\begin{itemize}
  \item Using Eqn.~(\ref{defn-avg}) and (\ref{defn-Q}), A is simply
    $\frac{\pd Q}{\pd t}$.

  \item For term B, we break up the sum and simplify each term seperately.
    The first term is
    \[
      \int^1_{-1} \int_0^{L_z} \int_0^{L_x} \frac{\pd}{\pd x}(U^2) \,d\bx
      =  \int^1_{-1} \int_0^{L_z}  U^2 \Bigg|_{x=0}^{x=L_x} 
      = 0
    \]
    since the flow is periodic in the x direction.  Similiarly, since the
    flow is also periodic in the z direction, we have that
    \[
      \int^1_{-1} \int_0^{L_z} \int_0^{L_x} \frac{\pd}{\pd z}(UW) \,d\bx
       = 0.
    \]
    For the third term, 
    \[
      \int^1_{-1} \int_0^{L_z} \int_0^{L_x} \frac{\pd}{\pd x}(UV) \,d\bx
        = \int^{L_x}_0 \int_0^{L_z}  UV \Bigg|_{y=-1}^{y=1} 
        = 0
    \]
    by the no--slip boundary condition along the walls. So, $B=0$.

  \item C must also be zero, for if the mean pressure gradient was not zero,
     then there would be no need to add an external force, $F_x$, to drive
     the flow.

  \item For D, we break up the sum as we did for B. So
    \[ 
      \int^1_{-1} \int_0^{L_z} \int_0^{L_x} \frac{\pd^2 U}{\pd x^2}
      = \int^1_{-1} \int_0^{L_z} \frac{\pd U}{\pd x} \Bigg|_{x=0}^{x=L_x}
      = 0
    \]
    since the flow is periodic in the x direction. A similiar argument can
    be used to show that
    \[ 
      \int^1_{-1} \int_0^{L_z} \int_0^{L_x} \frac{\pd^2 U}{\pd z^2}
      = 0.
    \]
    For the third term, we have that
    \begin{eqnarray*}
      \frac{L_z}{h}\frac{1}{L_xL_z} \int^1_{-1} \int_0^{L_z} \int_0^{L_x}
      \frac{1}{Re_0} \frac{\pd^2 U}{\pd y^2}
      &=& \frac{L_z}{h} \frac{1}{Re_0} 
        \int^1_{-1} \frac{\pd^2 <U>_{x,z}}{\pd y^2} \,dy
      \\
      &=& \frac{L_z}{h} \frac{1}{Re_0} 
        \frac{\pd <U>_{x,z}}{\pd y} \Bigg|_{y=-1}^{y=1}.
    \end{eqnarray*}
    Since the mean channel flow is symmetric about the channel center line,
    $y=0$,
    \[
      \frac{\pd <U>_{x,z}}{\pd y} \Bigg|_{y=1} 
      = -\frac{\pd <U>_{x,z}}{\pd y} \Bigg|_{y=-1}.
    \]
    So, 
    \[ 
      D = \frac{L_z}{h} \frac{1}{Re_0} 
                        \frac{\pd <U>_{x,z}}{\pd y} \Bigg|_{y=-1}^{y=1}
        = -2\frac{L_z}{h} \frac{1}{Re_0} 
                        \frac{\pd <U>_{x,z}}{\pd y} \Bigg|_{y=-1}.
    \]

  \item Finally, since $F_x$ is constant in space, E simplifies to
    \[
      \frac{L_z}{h}\frac{1}{L_xL_z} \int^1_{-1} \int_0^{L_z} \int_0^{L_x}
        F_x \,d\bx
      = \frac{2L_z}{h} F_x.
    \]

\end{itemize}
Using the above, Eqn.~(\ref{non-dim-moment-1}) can be rewritten as
\begin{equation}
  \label{Q-eqn}
  \frac{\pd Q}{\pd t} = -2\frac{L_z}{h} (u_\tau^2 - F_x),
\end{equation}
where $u_\tau^2 = \frac{1}{Re_0} \frac{\pd <U>_{x,z}}{\pd y} \Bigg|_{y=-1}$
is the mean wall shear stress.
So, as the flow changes from a laminar state to a turbulent state, $u_\tau^2$ 
increases and if the pressure gradient is kept constant at the laminar
value, Eqn.~(\ref{Q-eqn}) says that $Q$ will decrease. To maintain a constant
flow rate, $Q$, it is then necessary to balance the pressure gradient to 
the mean wall shear stress. However, to calculate the mean wall
shear stress would require a large grid to obtain the mean values in the 
stream and spanwise directions. This makes balancing the pressure gradient
with the mean wall shear stress computationally too expensive. So, $F_x$ is
found by relating the pressure gradient to the flow rate fluctuations from
the constant initial value $Q_0,$ 
\[
  F^{n+1} = F^n - a_1(Q^{n+1} - Q_0) - a_2(Q^n - Q_0)
\]
where $a_1 = 2\del t$, $a_2 = \del t$, and $n$ and $n+1$ represent the old
and new time states. Note that the formula given in the reference, \cite{tran},
contains typographical errors.

\subsection{Numerical Results}

Using the mean initial centerline velocity as our characteristic velocity and
the channel half height as our characteristic length scale, the Reynolds 
number for our numerical experiments is $Re^0_C=5000.$ Note that in the
following discussion that all mean calculations are based on averaging in 
the homogeneous directions (x and z) as well as time. Also, note that we use
as our characteristic velocity for non--dimensionalization, the wall 
shear velocity
\[
  u_\tau = \left( \frac{\tau_w}{\rho} \right)^\frac{1}{2}
\]
where $\tau_w$ is the wall shear stress and use as our non--dimensional
wall distance
\[
  y^+ = \frac{y u_\tau}{\nu}.
\]

\renewcommand{\textfraction}{0}

\begin{figure}[t]
\begin{center}
\leavevmode\epsfxsize=100mm
\leavevmode\epsfysize=80mm
\epsffile{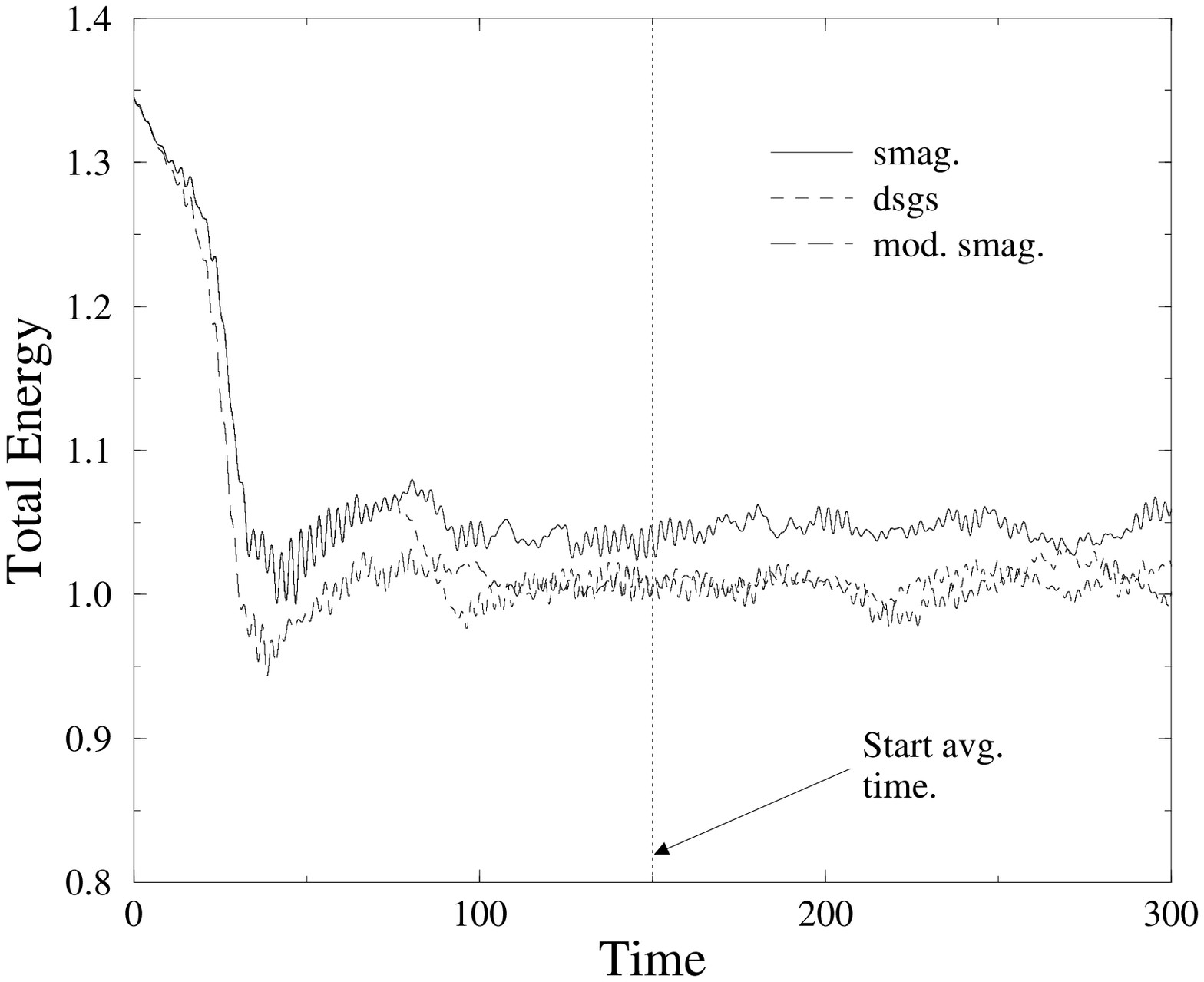}
\end{center}
\renewcommand{\baselinestretch}{1}
\caption{ \label{ke-fig} \em
          Plot of the total kinetic engery for the channel flow using the
          Smagorinsky, DSGS, and modified Smagorinsky subgrid scale models.
          At time $t = 150$, a running total of the fluid variables are started
          to calculate the various mean properties of the flow. }
\end{figure}

Fig.~\ref{ke-fig} is a plot of the total kinetic energy using the three
different models. As can be seen in the figure, by simulation time
$t=150$, the flow has stabilized enough to begin the time averaging
calculations. Note that in the DSGS model run, the Smagorinsky model is used
while the flow is developing until time $t=75.$ This is the reason for the
initial overlap of energy between the Smagorinsky and DSGS results in 
Fig.~\ref{ke-fig}.

\begin{figure}[h]
\begin{center}
\leavevmode\epsfxsize=100mm
\leavevmode\epsfysize=130mm
\epsffile{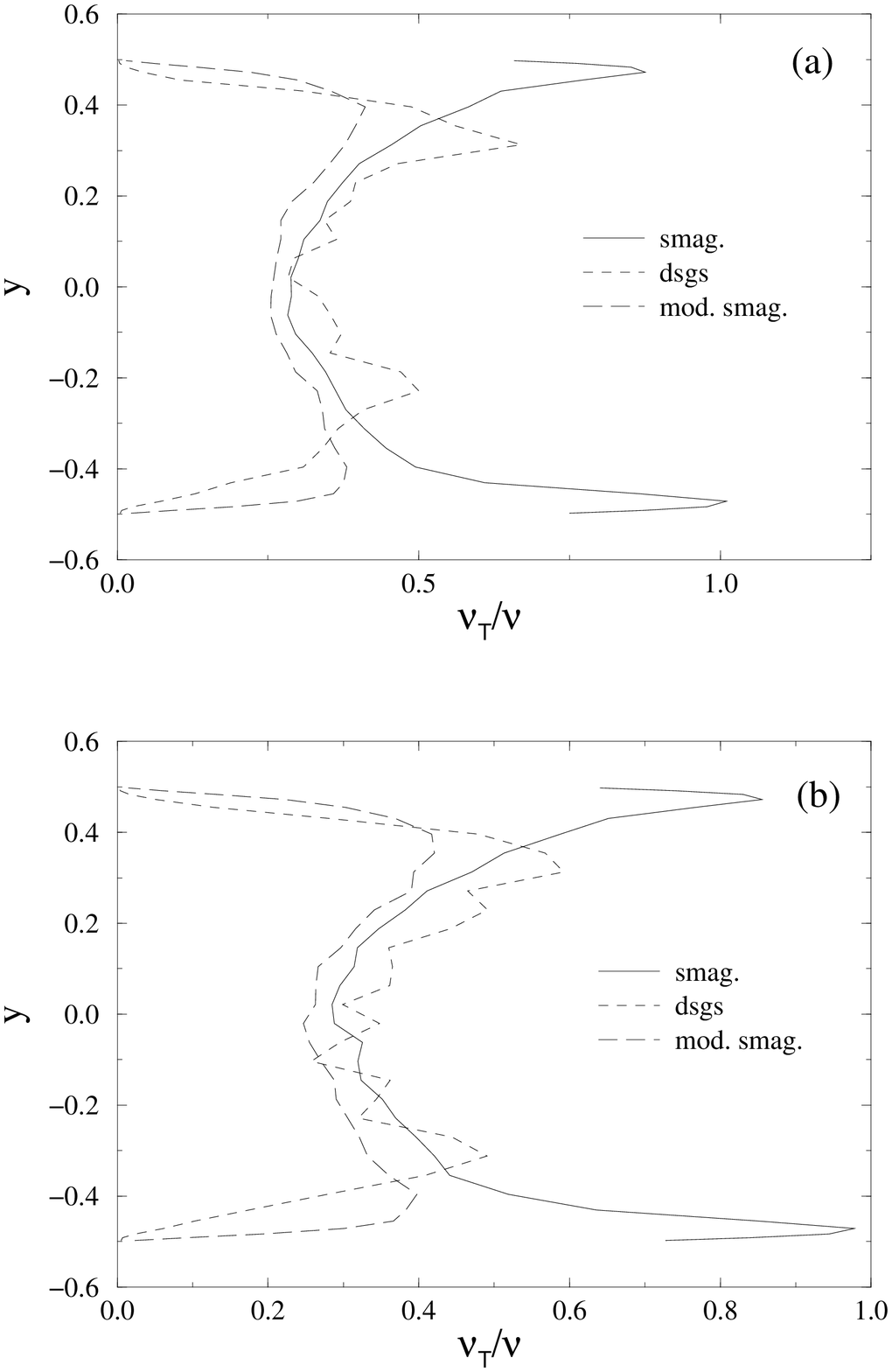}
\end{center}
\renewcommand{\baselinestretch}{1}
\caption{ \label{avg-nut-fig} \em
          Time averaged plots of the eddy viscosity normalized by the kinematic
          viscosity. (a) is a plot along the line $x = 5 h$ and $z = 1.0625 h.$
          (b) is a plot along the line $x = 9 h$ and $z = 1.0625 h$.}
\end{figure}
One of the problems with the Smagorinsky model discussed above is that the
eddy viscosity is not damped as the wall and the viscous sublayer is approached.
This can clearly be seen in Fig.~\ref{avg-nut-fig}, the sketch of the time
averaged eddy viscosity along the lines $x = 5 h$ and $z = 1.0625 h$, and
$x = 9 h$ and $z = 1.0625 h$. 
\begin{figure}[h]
\begin{center}
\leavevmode\epsfxsize=100mm
\leavevmode\epsfysize=130mm
\epsffile{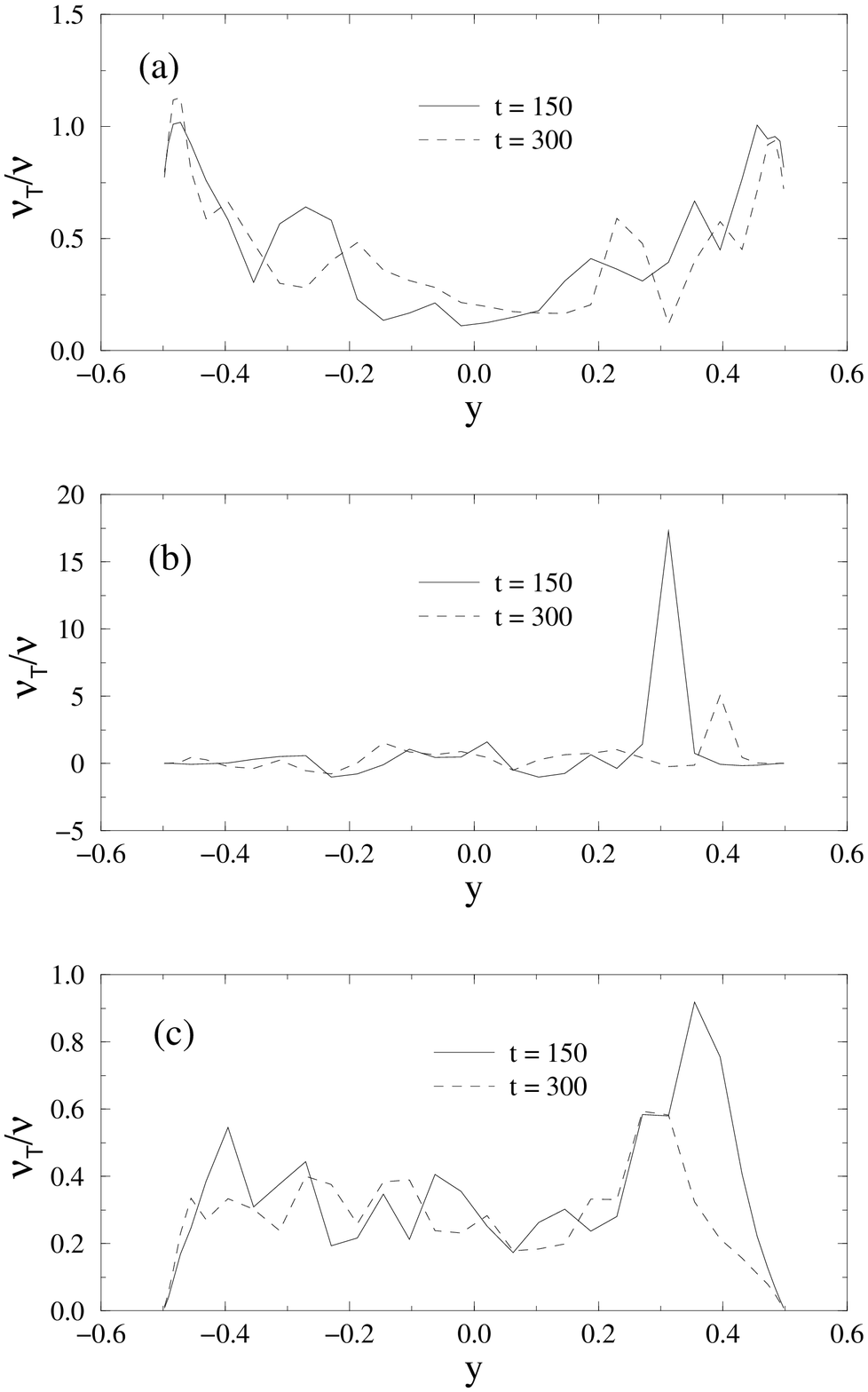}
\end{center}
\renewcommand{\baselinestretch}{1}
\caption{ \label{nut-fig} \em
          Instantaneous snap shot of the eddy viscosity normalized by the 
          kinematic viscosity along the line $x = 5 h$ and $z = 1.0625 h.$
          (a) is using the Smagorinsky model, (b) is the DSGS model, and
          (c) is the modified Smagorinsky model. }
\end{figure}
Note that along these lines, the magnatude of 
the eddy viscosity is approximately 75\% of the kinematic viscosity at the 
walls, when using the Smagorinsky model. The plots also show that both
the DSGS and the modified Smagorinsky model damp out the eddy viscosity as
the wall is approached. Fig.~\ref{avg-nut-fig} also indicates that the
damping function in the modified Smagorinsky model not only damps out the
eddy viscosity near the walls but also causes some damping far from the walls.
However, the plots show that the damping along the central region of the
channel is not as great as the damping near the walls. Fig.~\ref{nut-fig},
a snap shot of the eddy viscosity along the line $x = 5 h$ and 
$z = 1.0625 h$ at simulation times $t = 150$ and $t = 300$, shows that at any
given time, the eddy viscosity can be over a magnatude greater than the 
kinematic viscosity when using the DSGS model. Furthermore, Fig.~\ref{nut-fig}b
shows that at any given time, the eddy viscosity can be negative (which can
be interpreted as backscatter) when using the DSGS model.

\begin{figure}[h]
\begin{center}
\leavevmode\epsfxsize=80mm
\leavevmode\epsfysize=60mm
\epsffile{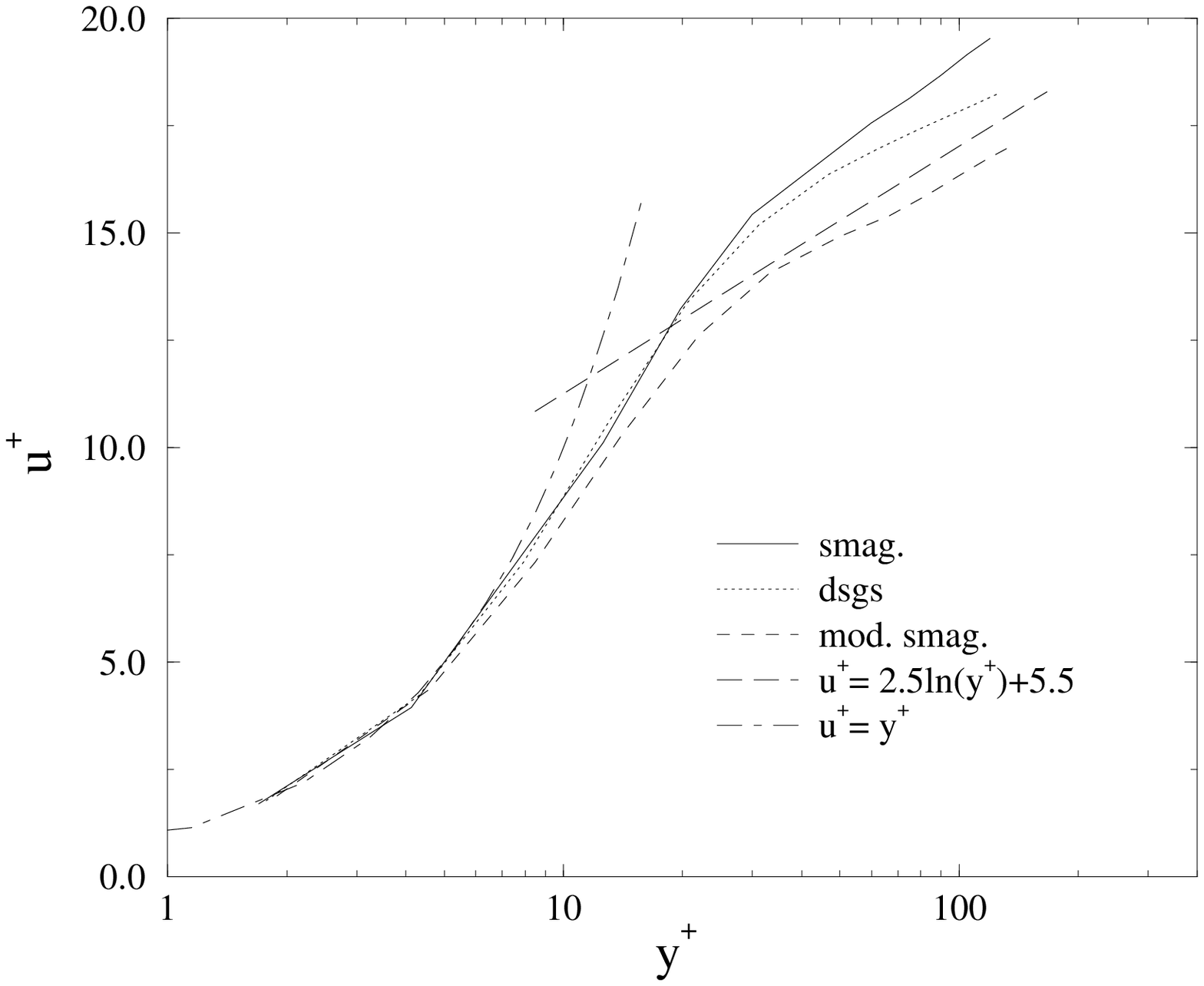}
\end{center}
\renewcommand{\baselinestretch}{1}
\caption{ \label{u-mean-fig} \em
          Log--linear plot of the mean velocity normalized by the wall shear
          velocity, $u_\tau,$ for the 
          channel flow using the Smagorinsky, DSGS, and modified Smagorinsky 
          subgrid scale models. }
\end{figure}

Fig.~\ref{u-mean-fig}, a log--linear plot of the non--dimensionalized mean 
velocity, indicates that the modified Smagorinsky data is much closer
to the logarithmic friction law than the Smagorinsky model. Note that we
use as our additive constant in the log law, $5.5$, as in Kim, Moin, and Moser's
work \cite{kim} since our Reynolds number, based on the wall shear velocity,
is $Re_\tau \approx 180$ for our channel computations. 

\begin{figure}[h]
\begin{center}
\leavevmode\epsfxsize=80mm
\leavevmode\epsfysize=60mm
\epsffile{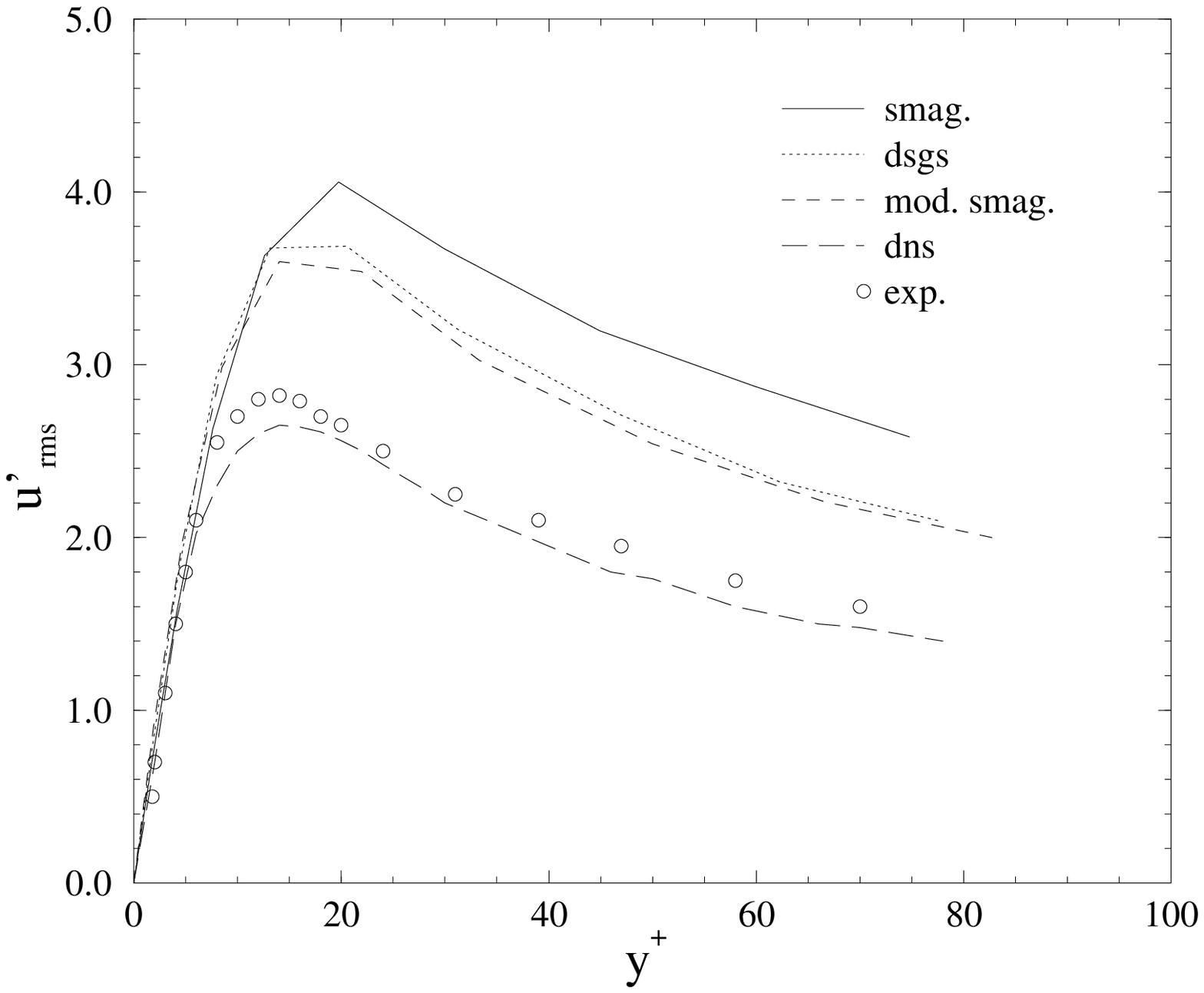}
\end{center}
\renewcommand{\baselinestretch}{1}
\caption{ \label{u-rms-fig} \em
          Root-mean-square of the streamwise velocity fluctuations normalized
          by the wall shear velocity. Note that dns data is from 
          Kim~et~al.~\cite{kim} and the experimental data, exp, is from
          Kreplin and Eckelmann \cite{kreplin}. }
\end{figure}

\begin{figure}[h]
\begin{center}
\leavevmode\epsfxsize=80mm
\leavevmode\epsfysize=60mm
\epsffile{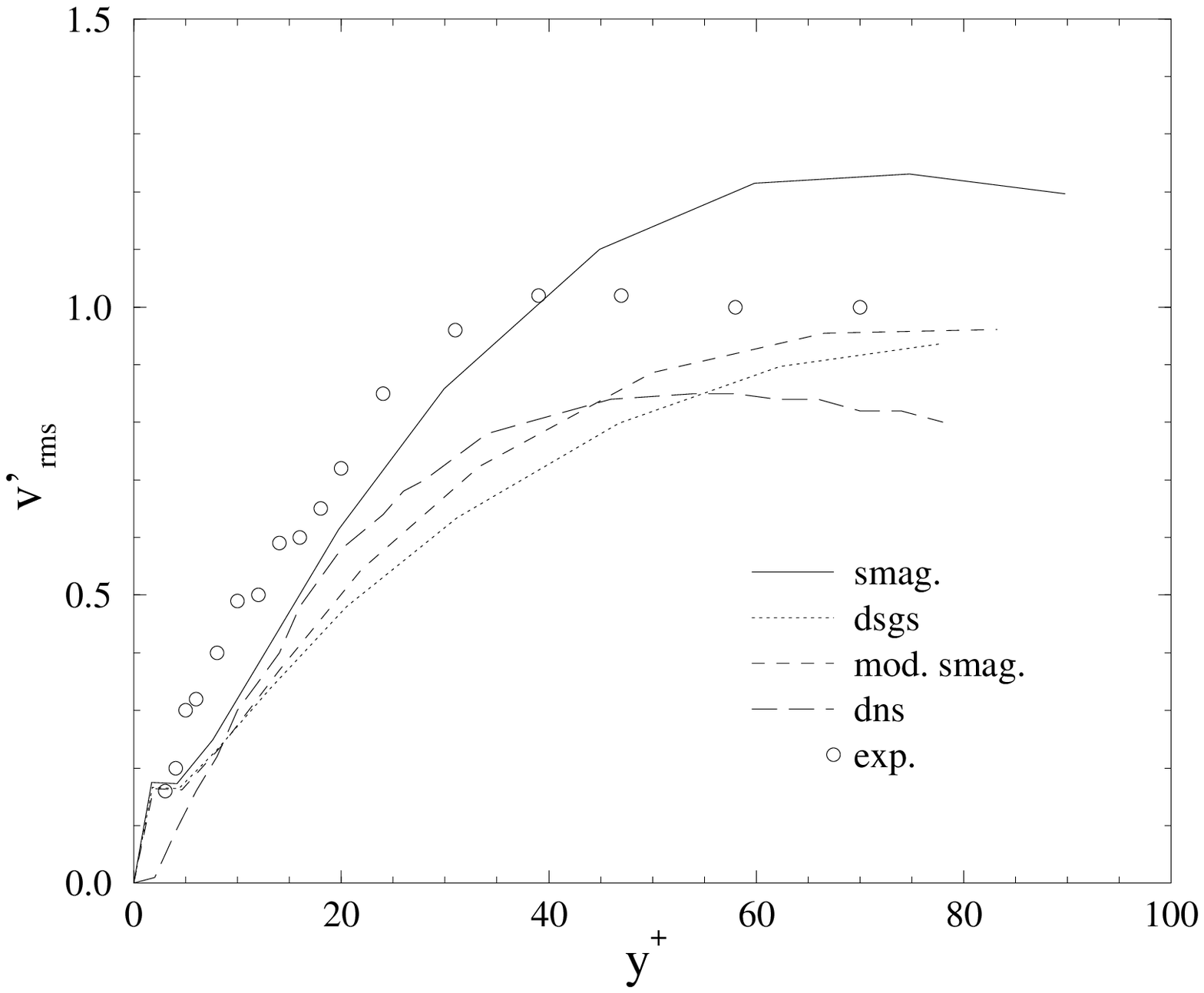}
\end{center}
\renewcommand{\baselinestretch}{1}
\caption{ \label{v-rms-fig} \em
          Root-mean-square of the wall normal velocity fluctuations normalized
          by the wall shear velocity. Note that dns data is from 
          Kim~et~al.~\cite{kim} and the experimental data, exp, is from
          Kreplin and Eckelmann \cite{kreplin}. }
\end{figure}

\begin{figure}[h]
\begin{center}
\leavevmode\epsfxsize=80mm
\leavevmode\epsfysize=60mm
\epsffile{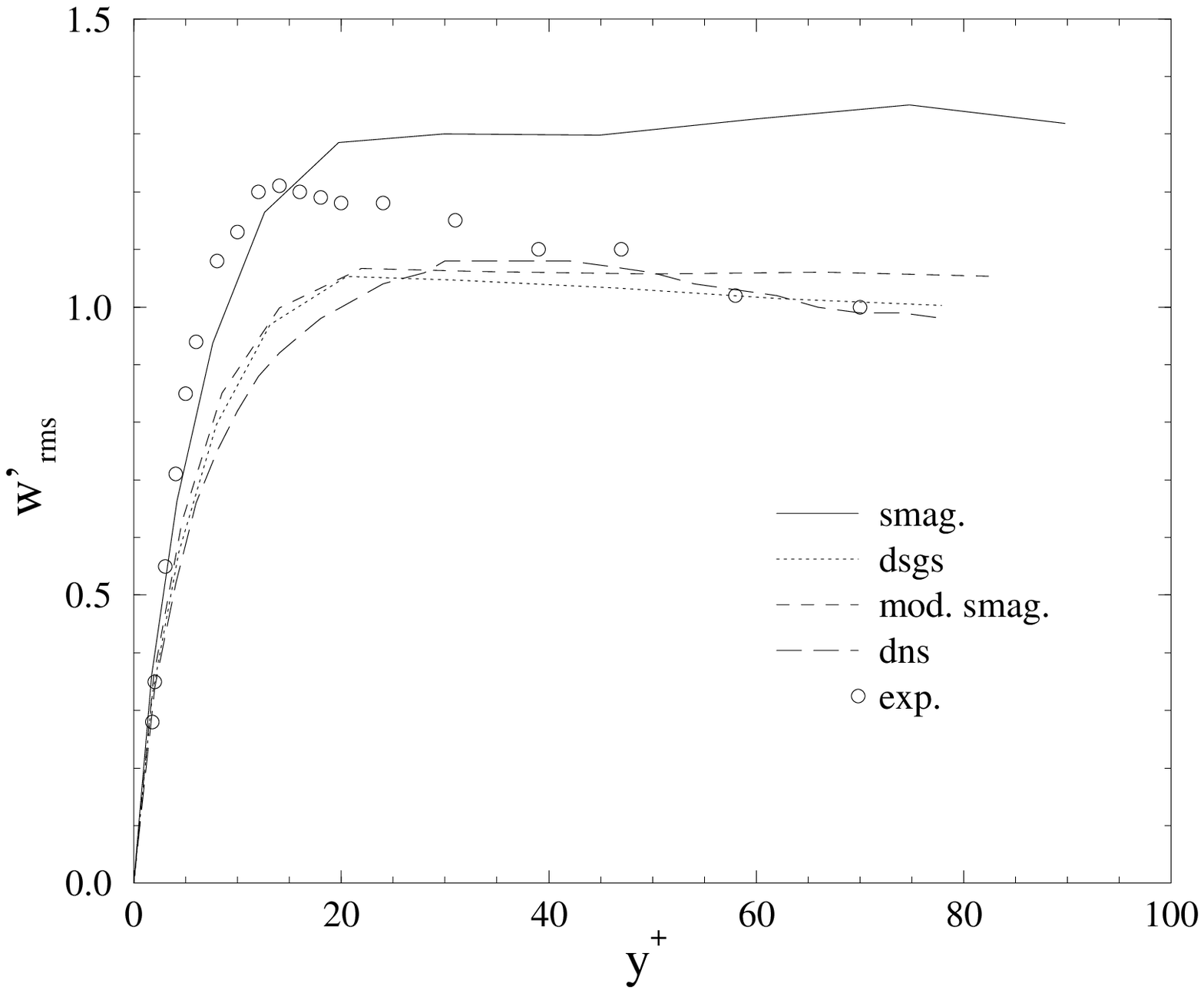}
\end{center}
\renewcommand{\baselinestretch}{1}
\caption{ \label{w-rms-fig} \em
          Root-mean-square of the spanmwise velocity fluctuations normalized
          by the wall shear velocity. Note that dns data is from 
          Kim~et~al.~\cite{kim} and the experimental data, exp, is from
          Kreplin and Eckelmann \cite{kreplin}. }
\end{figure}

Figures \ref{u-rms-fig} -- \ref{w-rms-fig} are plots of the root-mean-square
of the velocity fluctuations. Note that the DNS results to which we compare
our data is that of Kim~et~al.~\cite{kim} and the experimental data is that
of Kreplin and Eckelmann~\cite{kreplin}. It is not very surprising that the
velocity fluctuation properties of our LES computations do not match the 
DNS and experimental results exactly, because of the way in which the 
fluctuations are computed. If we let $<\cdot>_{t,x,z}$ be the average in time
and the homogeneous directions and $'$ denotes the statistical fluctuations
with zero mean, then note that
\[ 
  (\ua')^2 = <\ua>_{t,x,z}^2 - <\ua^2>.
\]
However, for our LES cases, the only values we have are at the filter level,
$\fua.$ Hence the only fluctuation we can calculate is
\[
  (\fua')^2 = <\fua>_{t,x,z}^2 - <\fua^2>_{t,x,z}
\]
where $\fua'$ represents the statistical fluctuations in the filtered field.
Note that $\fua'$ is well defined (in the sense that it will not be identically
zero) since, as discussed in the filter section, 
only the Sharp--Cutoff filter has compact support in wave space and so for
filters such as the Top--Hat filter, $\fua$ will still contain some small 
scale motions. Also, the grading of the grid near the walls implies that the
filter width decreases near the walls and so will contain more of the small
scale motions. Hence, we would expect $\fua'$ to be closer to $\ua'$ near the
walls, were the filter width is much smaller than toward the channel center.
Figures \ref{u-rms-fig} -- \ref{w-rms-fig} show that this is indeed the case. 

\begin{figure}[h]
\begin{center}
\leavevmode\epsfxsize=80mm
\leavevmode\epsfysize=60mm
\epsffile{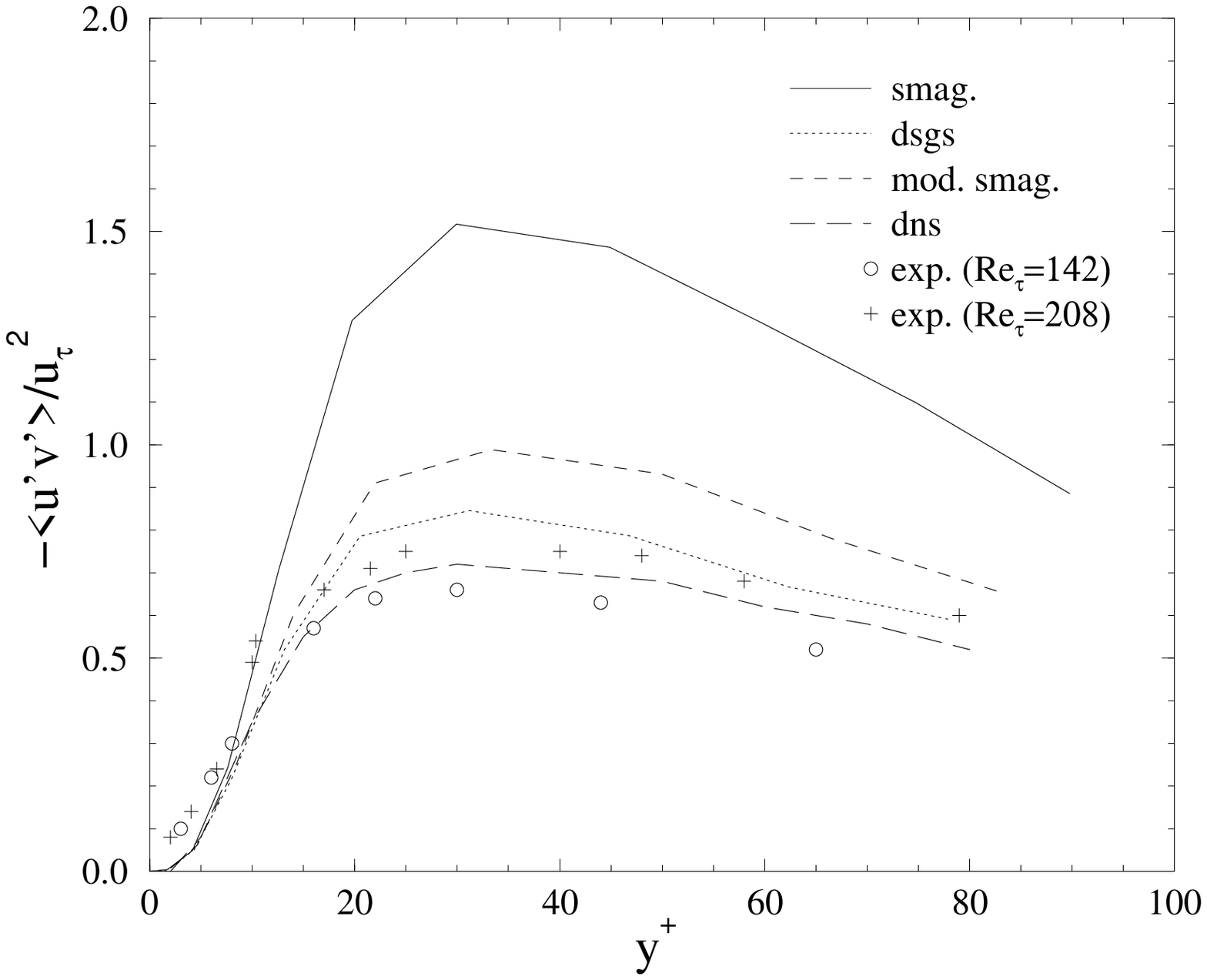}
\end{center}
\renewcommand{\baselinestretch}{1}
\caption{ \label{uv-fig} \em
          Reynolds shear stress normalized by the wall shear velocity.
          Note that dns data is from 
          Kim~et~al.~\cite{kim} and the experimental data, exp, is from
          Eckelmann \cite{eckelmann} at two different Reynolds numbers, 
          $Re_\tau = 142$ and $Re_\tau = 208$. }
\end{figure}

Fig.~\ref{uv-fig} is a plot of the Reynolds shear stress normalized by the 
wall shear velocity. As discussed above, we note that the LES results are 
much closer to the DNS results near the wall. Also note that the DSGS and
modified Smagorinsky are closer to the DNS results than using standard
Smagorinsky model.

\begin{figure}[h]
\begin{center}
\leavevmode\epsfxsize=80mm
\leavevmode\epsfysize=60mm
\epsffile{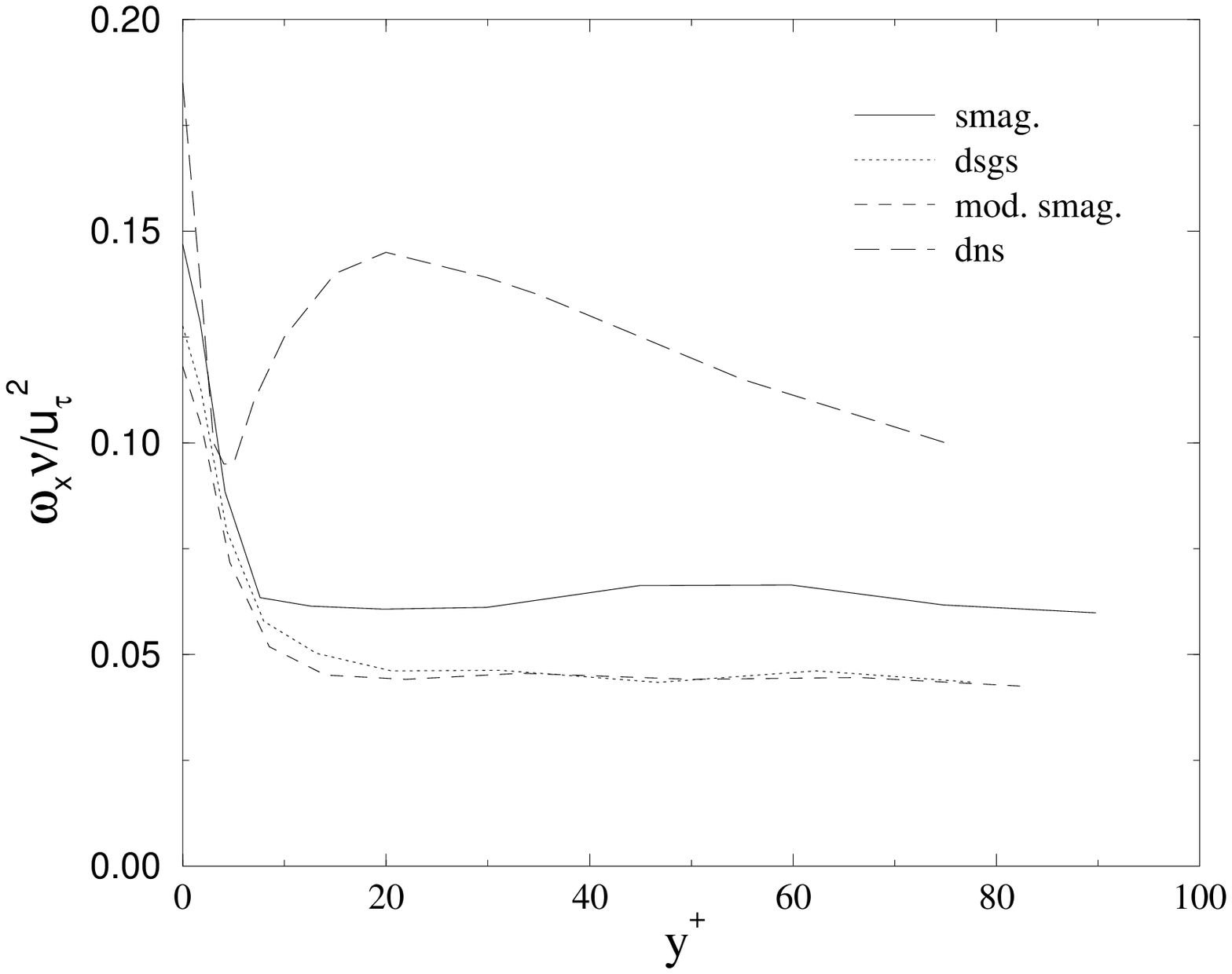}
\end{center}
\renewcommand{\baselinestretch}{1}
\caption{ \label{vort-x-fig} \em
          Root-mean-square vorticity (x) fluctuations normalized
          by the mean wall shear. Note that dns data is from 
          Kim~et~al.~\cite{kim}. }
\end{figure}

\begin{figure}[h]
\begin{center}
\leavevmode\epsfxsize=80mm
\leavevmode\epsfysize=60mm
\epsffile{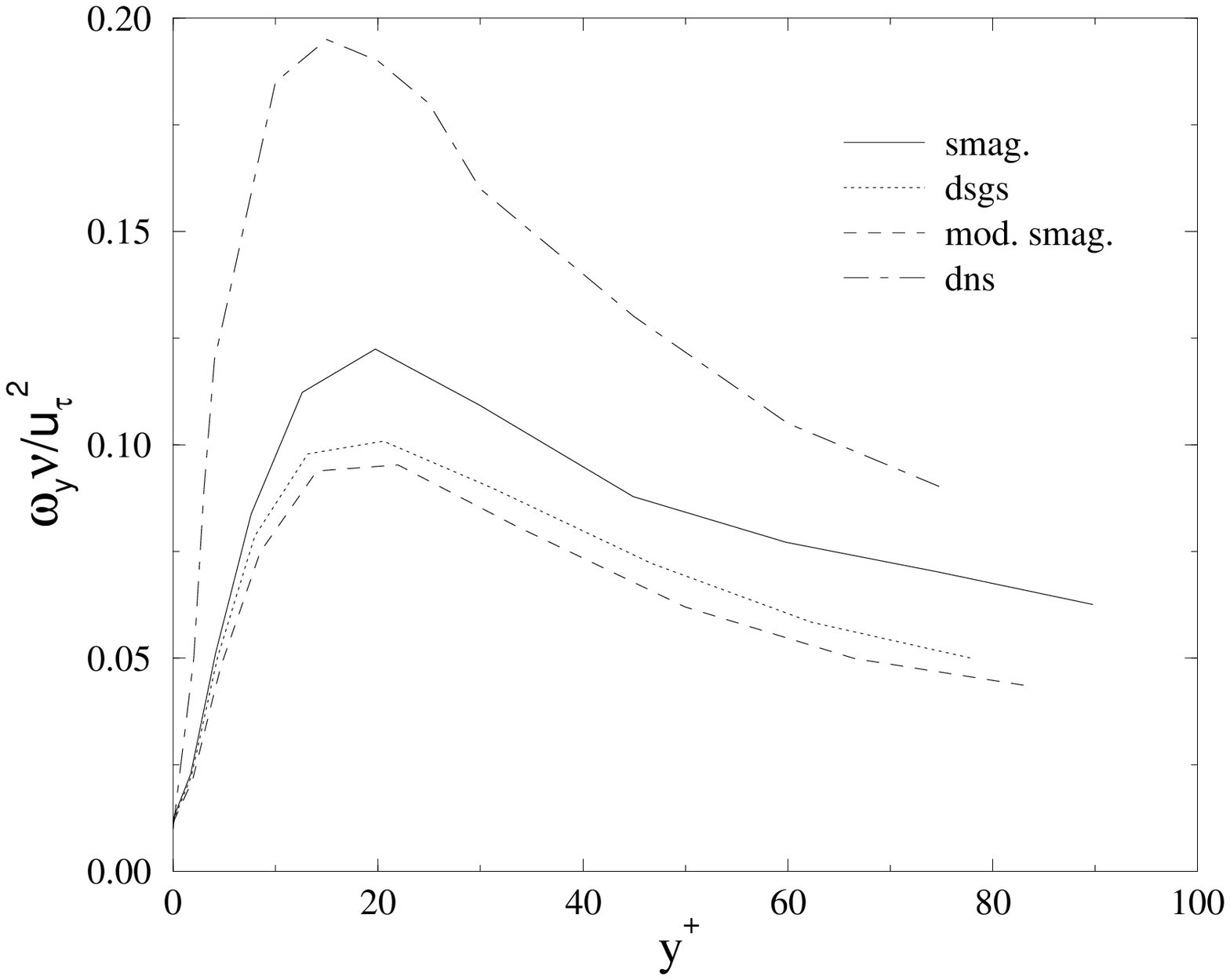}
\end{center}
\renewcommand{\baselinestretch}{1}
\caption{ \label{vort-y-fig} \em
          Root-mean-square vorticity (y) fluctuations normalized
          by the mean wall shear. Note that dns data is from 
          Kim~et~al.~\cite{kim}. }
\end{figure}

\begin{figure}[h]
\begin{center}
\leavevmode\epsfxsize=80mm
\leavevmode\epsfysize=60mm
\epsffile{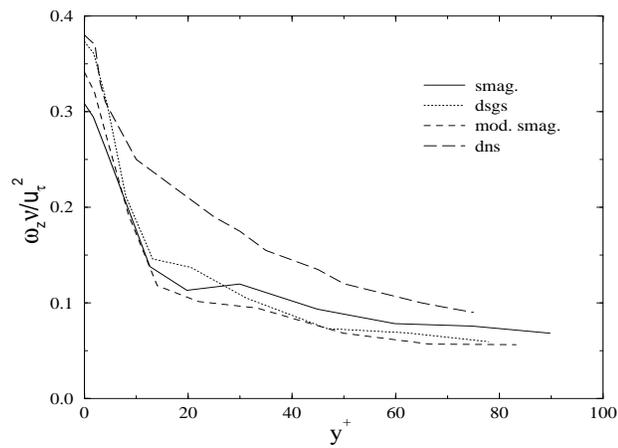}
\end{center}
\renewcommand{\baselinestretch}{1}
\caption{ \label{vort-z-fig} \em
          Root-mean-square vorticity (z) fluctuations normalized
          by the mean wall shear. Note that dns data is from 
          Kim~et~al.~\cite{kim}. }
\end{figure}

Figures \ref{vort-x-fig} -- \ref{vort-z-fig} are plots of the root-mean-square
vorticity fluctuations normalized by the mean wall shear. Since vorticity
is mostly associated with small scale motions, it is not surprising that
the LES vorticity data are all lower than the DNS data for the reason given
above. Note that near the wall, at least one of the LES data is close to
the DNS results. Unfortunately, there is no consistency as to which model
is closest to the DNS results. For x vorticity fluctuations, the Smagorinsky
model is closest to the DNS data near the wall, the modified Smagorinsky
is the furthest off from the DNS data, and the DSGS model is between the
other two. However, for the z vorticity, note that the order of performance
of the models changes to the DSGS beging the closest, followed by the modified
Smagorinsky and then the Smagorinsky. 

\clearpage
\newpage
\section{Flow over a Backward Facing Step}

To study the effectiveness of the various SGS models on turbulent flows
with separation and reattachment, the flow over a backward facing step
is used because of its geometric simplicity. In our modeling of this 
flow, we use periodic boundary conditions in the spanwise direction. 

\subsection{Computational Parameters}

\begin{figure}[t]
\begin{center}
\leavevmode\epsfxsize=5in
\leavevmode\epsfysize=3in
\epsffile{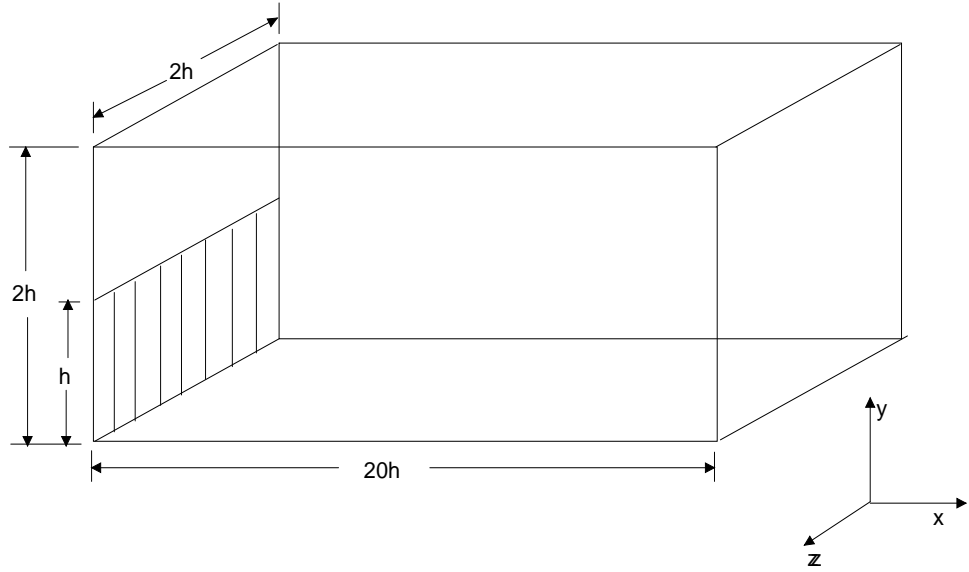}
\end{center}
\renewcommand{\baselinestretch}{1}
\caption{ \label{bfs-sketch} \em
          Sketch of the backward facing step and channel with approximate 
          relevant dimensions. For our experiment, the step height 
          is $h=0.45833$}
\end{figure}

The computational domain of our backward facing step problem consists of one
main channel (post expansion section) with the step located at one end of 
the channel, 
see Fig.~\ref{bfs-sketch}. The discretization of the channel consists of
$88 \times 32 \times 16$ elements and the dimensions of the channel are
approximately 20h in the streamwise direction and 2h in the spanwise and 
wall normal directions, where $h = 0.45833$ is the step height.

At the inlet, the fluid is forced by setting the streamwise velocity to be
1 along the inlet plane. Thus the Reynolds number, based on the inlet
velocity and the step height, is $Re_h=4583$. 

\subsection{Numerical Results}

\begin{figure}[b]
\begin{center}
\leavevmode\epsfxsize=6in
\leavevmode\epsfysize=6in
\epsffile{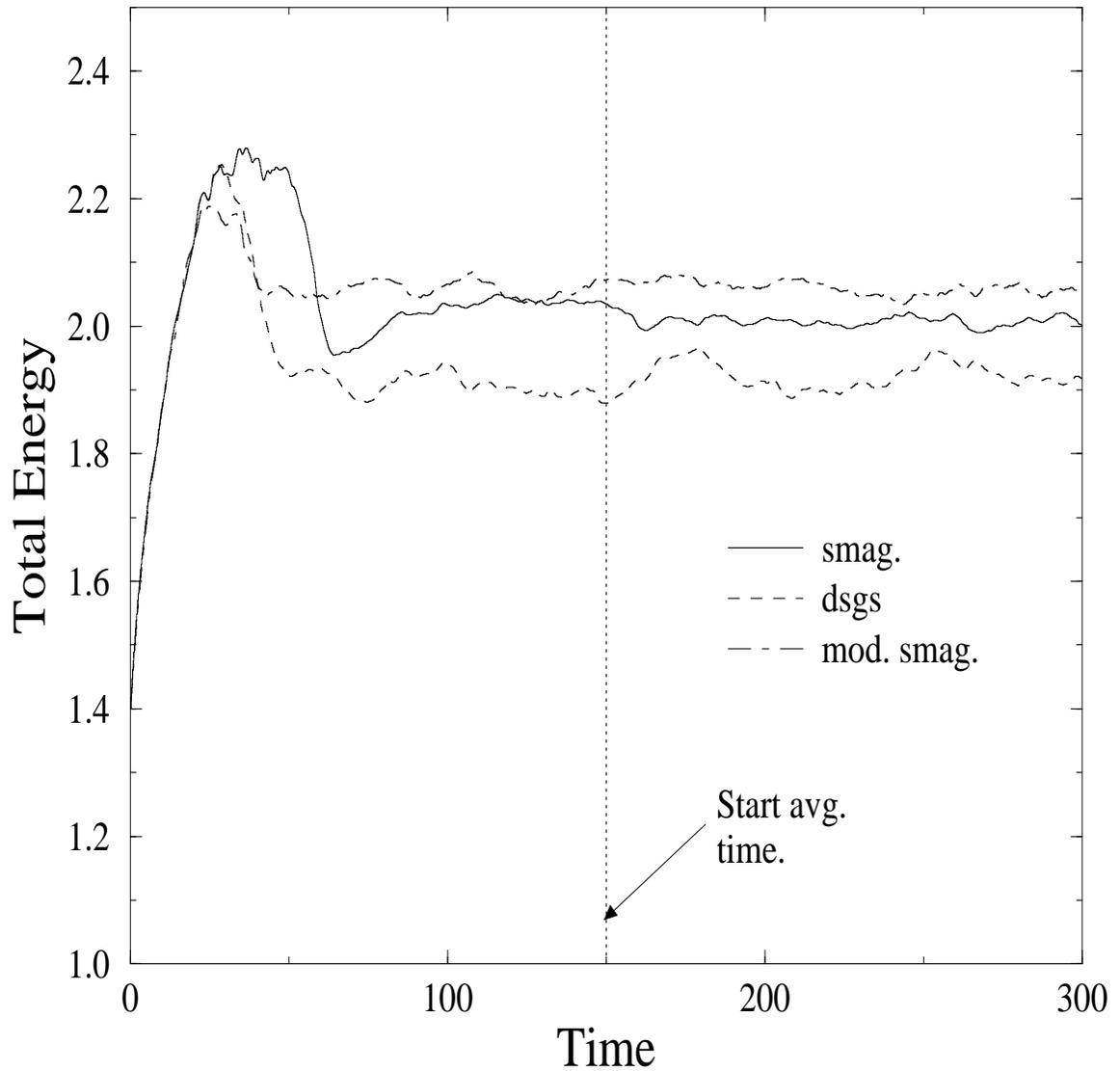}
\end{center}
\renewcommand{\baselinestretch}{1}
\caption{ \label{bfs-ke} \em
          Plot of the total kinetic energy for the backward facing step using
          the Smagorinsky, DSGS, and modified Smagorinsky subgrid scale 
          models. At time $t = 150$, a running total of the fluid variables
          is started to calculate the various mean properties of the flow. }
\end{figure}

\begin{figure}[b]
\begin{center}
\leavevmode\epsfxsize=6in
\leavevmode\epsfysize=6in
\epsffile{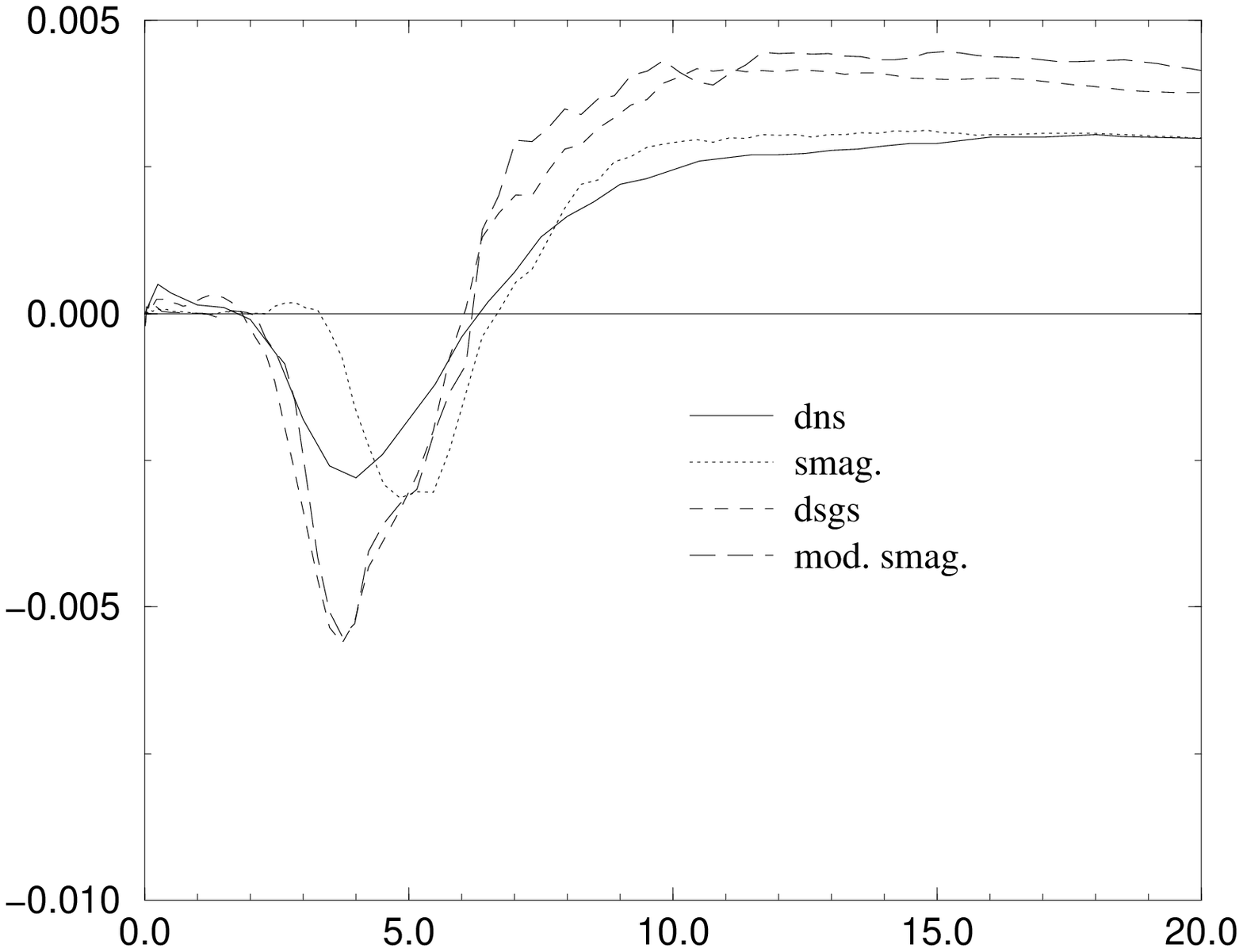}
\end{center}
\renewcommand{\baselinestretch}{1}
\caption{ \label{bfs-skin} \em
          Skin friction coefficient along the bottom wall of the post 
          expansion section of the backward facing step problem. dns:
          direct numercial simulation data of Le and Moin \cite{le},
          smag.: original Smagorinsky subgrid scale model, dsgs: dynamic
          subgrid scale model, and mod.~smag.: modified Smagorinsky model.}
\end{figure}

As in the channel flow simulations, the time averaging calculations are
started once the total kinetic energy has reached a quasi--equilibrium state, 
see Fig.~\ref{bfs-ke}. Figure \ref{bfs-skin} is a plot of the mean skin friction
coefficient,
\[ C_f = \frac{\tau_w}{\frac{1}{2}\rho U_0^2} \]
along the bottom wall, where $U_0$ is the inlet velocity. According to 
Le \& Moin \cite{le}, the high $|C_f|$ in the recirculation region is the
result of not having an entry section. 

\begin{figure}[b]
\begin{center}
\leavevmode\epsfxsize=6in
\leavevmode\epsfysize=7.5in
\epsffile{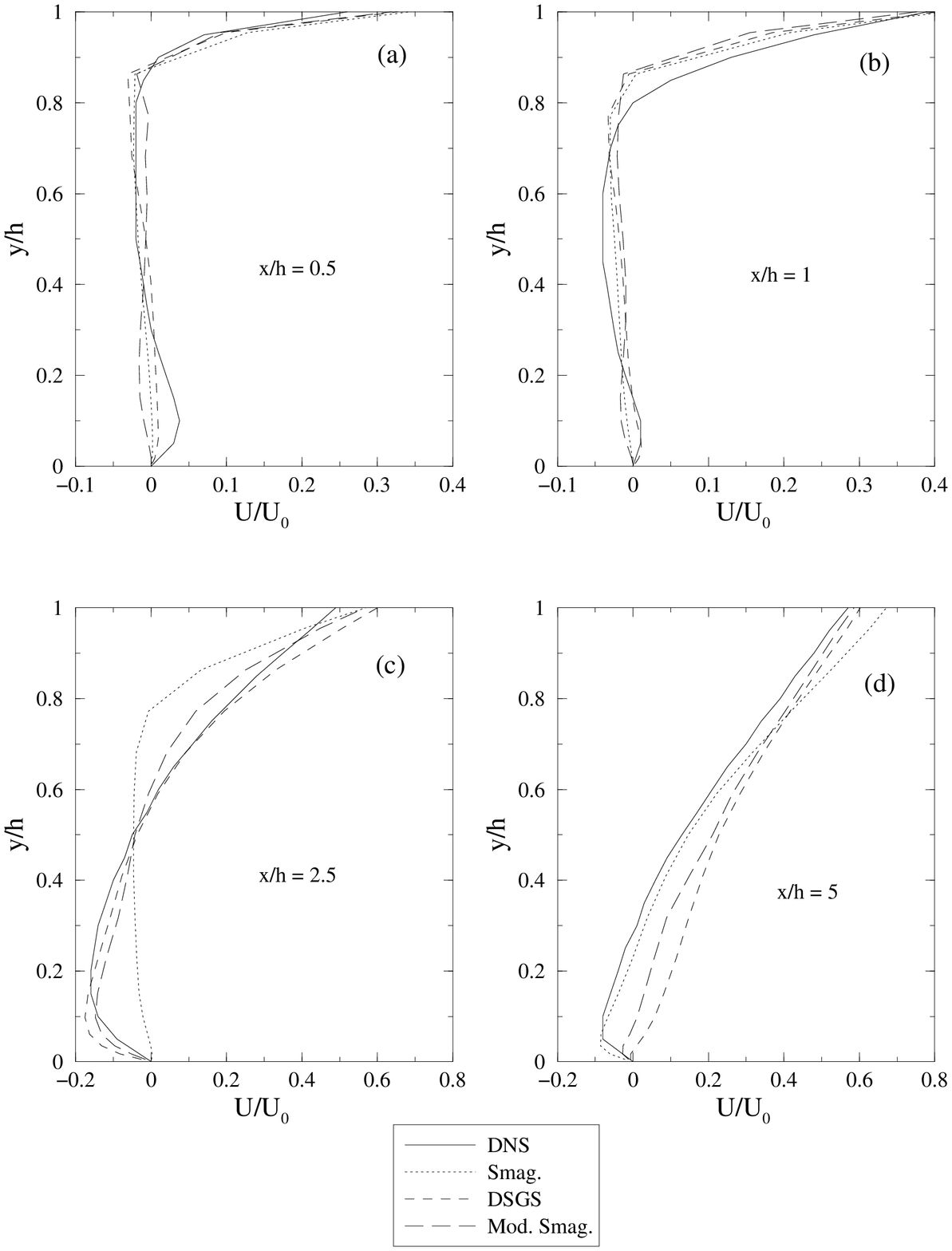}
\end{center}
\renewcommand{\baselinestretch}{1}
\caption{ \label{bfs-u1} \em
          Mean streamwise velocity profiles in the post 
          expansion section of the backward facing step problem. dns:
          direct numercial simulation data of Le and Moin \cite{le},
          smag.: original Smagorinsky subgrid scale model, dsgs: dynamic
          subgrid scale model, and mod.~smag.: modified Smagorinsky model.}
\end{figure}

\begin{figure}[h]
\begin{center}
\leavevmode\epsfxsize=6in
\leavevmode\epsfysize=7.5in
\epsffile{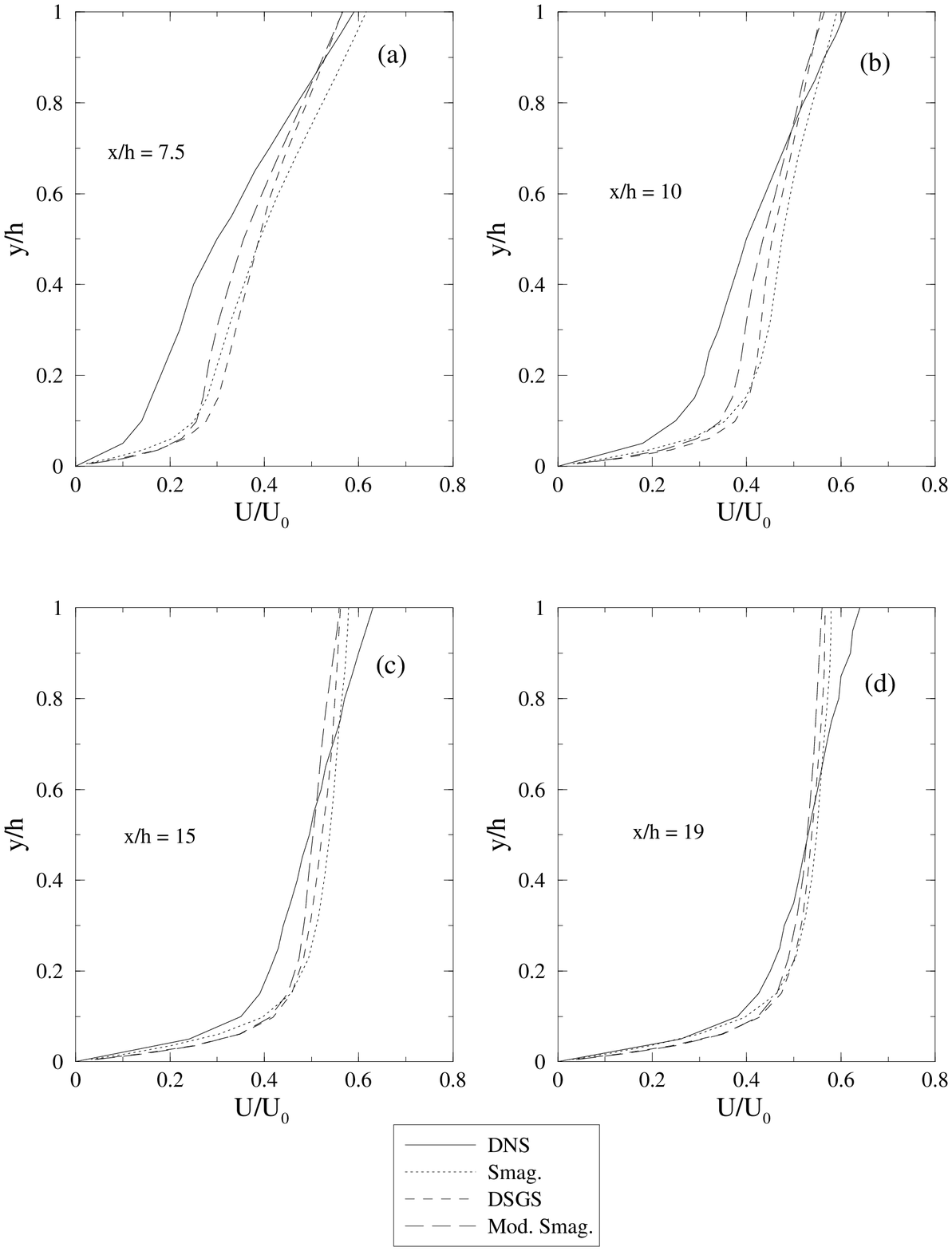}
\end{center}
\renewcommand{\baselinestretch}{1}
\caption{ \label{bfs-u2} \em
          Mean streamwise velocity profiles in the post 
          expansion section of the backward facing step problem. dns:
          direct numercial simulation data of Le and Moin \cite{le},
          smag.: original Smagorinsky subgrid scale model, dsgs: dynamic
          subgrid scale model, and mod.~smag.: modified Smagorinsky model.}
\end{figure}

Figures \ref{bfs-u1} and \ref{bfs-u2} are comparisons of the mean streamwise
velocity at various locations in the recirculation, reattachment, and 
recovery regions. In the reattachment and recovery regions, the data 
indicate that near the wall all the LES models overpredicted the mean streamwise 
velocity. In the region of the secondary bubble, the skin friction data
(Fig.~\ref{bfs-skin}) and the mean velocity data (Fig.~\ref{bfs-u1}) indicate that for the Smagorinsky
model and the modified Smagorinsky model, the secondary bubble is not as
well defined as in the DNS and DSGS cases.

\begin{figure}[h]
\begin{center}
\leavevmode\epsfxsize=5in
\leavevmode\epsfysize=6in
\epsffile{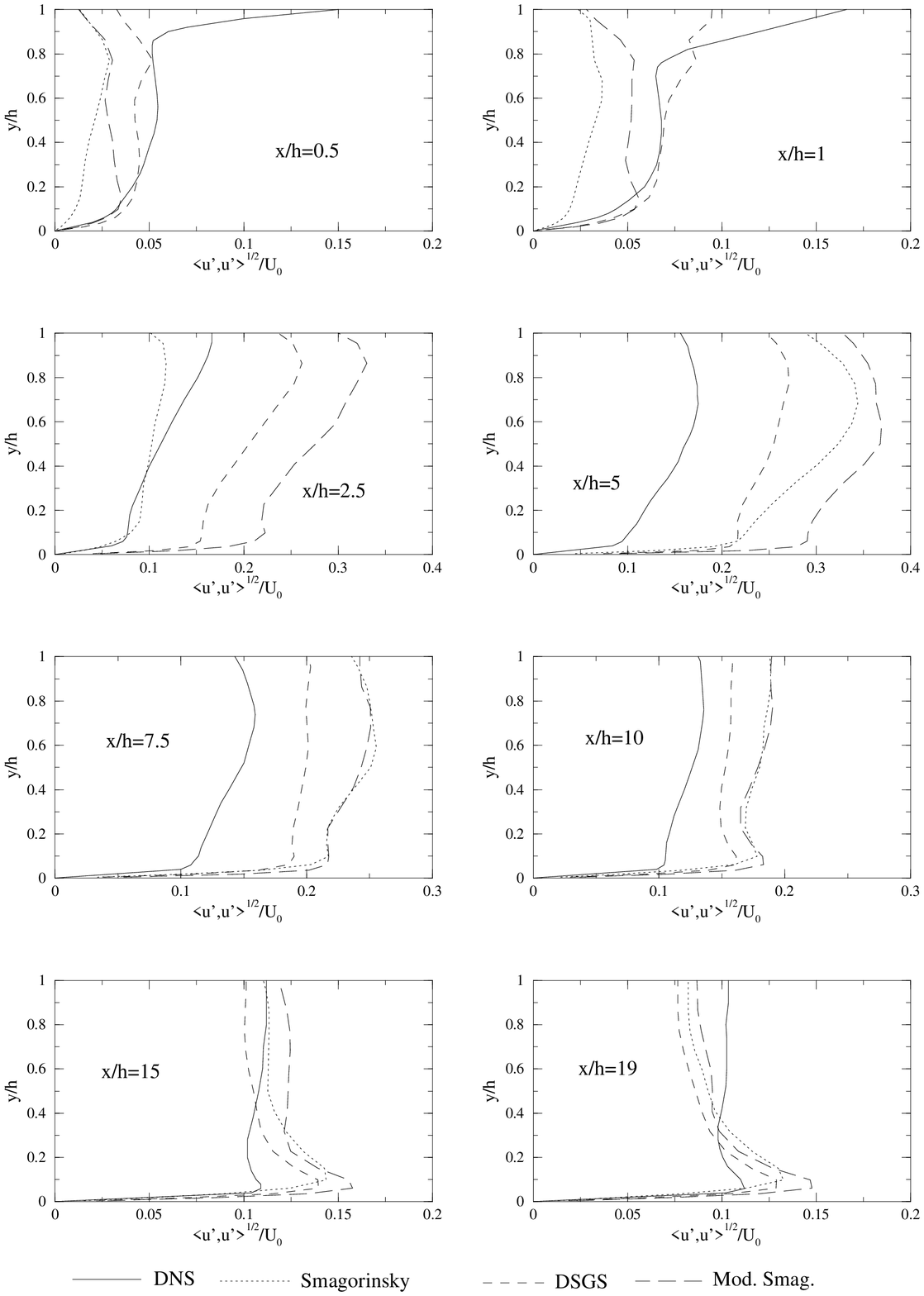}
\end{center}
\renewcommand{\baselinestretch}{1}
\caption{ \label{bfs-u-rms} \em
          Streamwise turbulence intensity profiles in the post 
          expansion section of the backward facing step problem. DNS:
          direct numercial simulation data of Le and Moin \cite{le},
          Smagorinsky: original Smagorinsky subgrid scale model, DSGS: dynamic
          subgrid scale model, and Mod.~Smagorinsky: modified Smagorinsky model.}
\end{figure}

\begin{figure}[h]
\begin{center}
\leavevmode\epsfxsize=5in
\leavevmode\epsfysize=6.5in
\epsffile{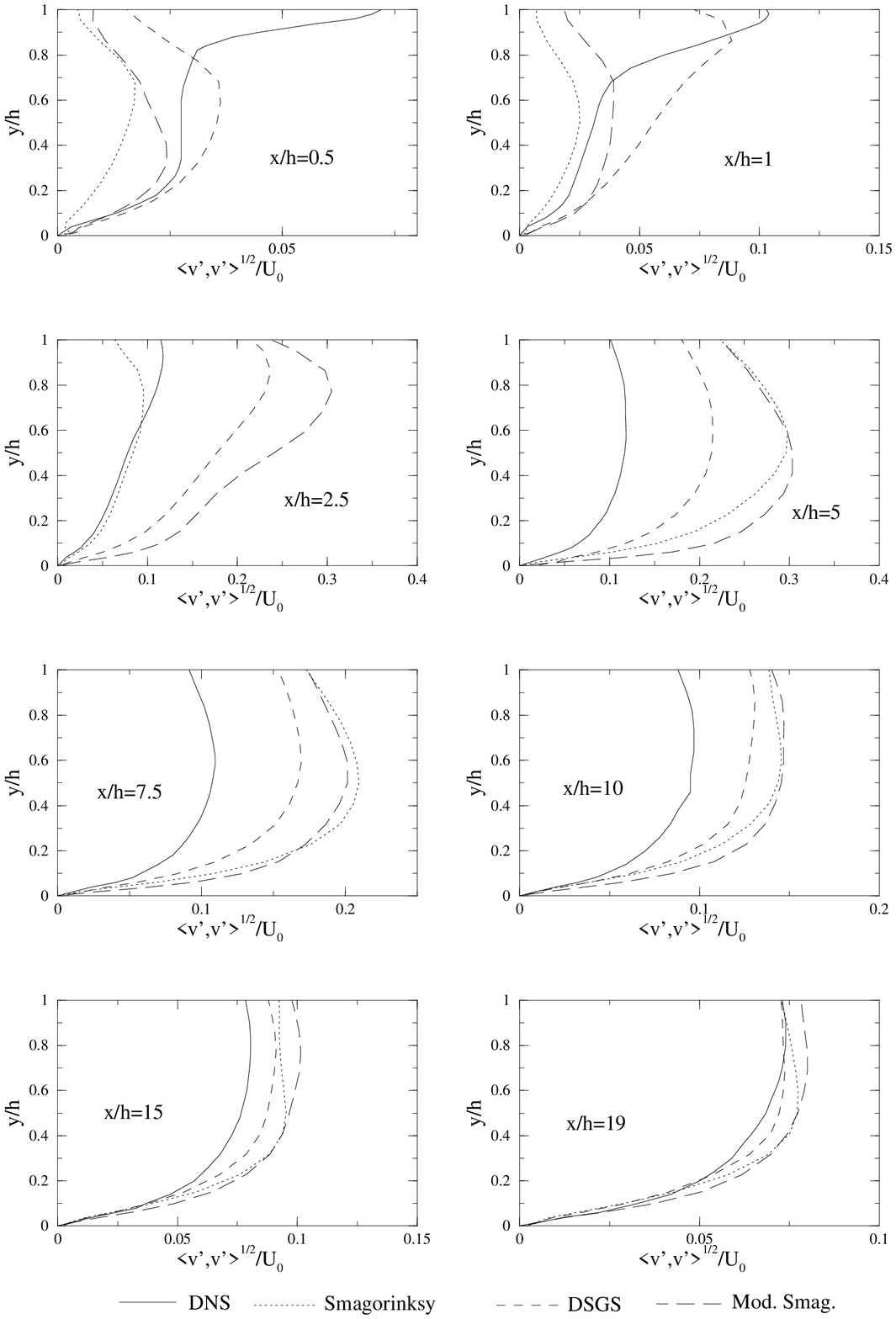}
\end{center}
\renewcommand{\baselinestretch}{1}
\caption{ \label{bfs-v-rms} \em
          Wall normal turbulence intensity profiles in the post 
          expansion section of the backward facing step problem. DNS:
          direct numercial simulation data of Le and Moin \cite{le},
          Smagorinsky: original Smagorinsky subgrid scale model, DSGS: dynamic
          subgrid scale model, and Mod.~Smagorinsky: modified Smagorinsky model.}
\end{figure}

\begin{figure}[h]
\begin{center}
\leavevmode\epsfxsize=5in
\leavevmode\epsfysize=6in
\epsffile{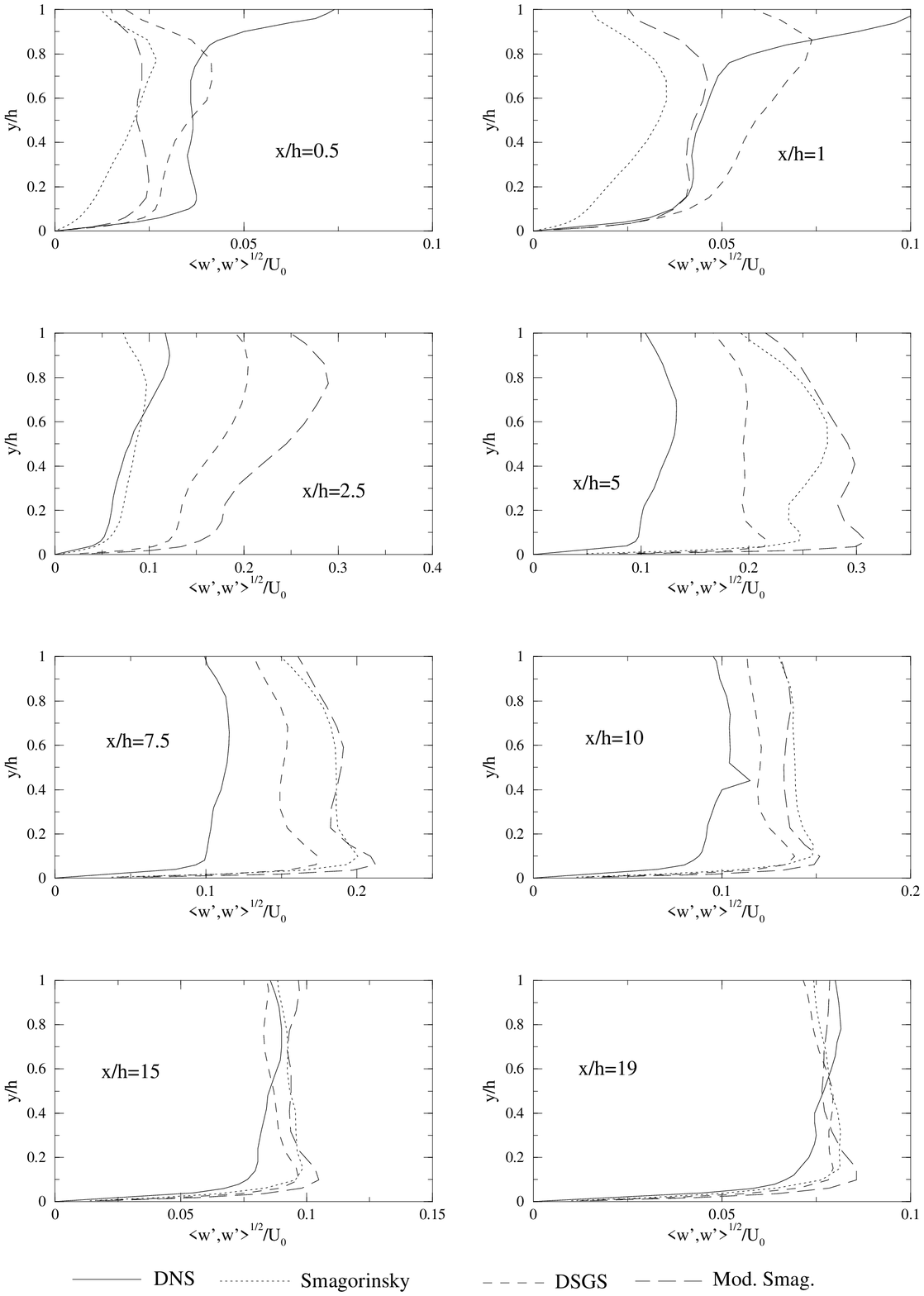}
\end{center}
\renewcommand{\baselinestretch}{1}
\caption{ \label{bfs-w-rms} \em
          Spanwise turbulence intensity profiles in the post 
          expansion section of the backward facing step problem. DNS:
          direct numercial simulation data of Le and Moin \cite{le},
          Smagorinsky: original Smagorinsky subgrid scale model, DSGS: dynamic
          subgrid scale model, and Mod.~Smagorinsky: modified Smagorinsky model.}
\end{figure}

\begin{figure}[h]
\begin{center}
\leavevmode\epsfxsize=5in
\leavevmode\epsfysize=6in
\epsffile{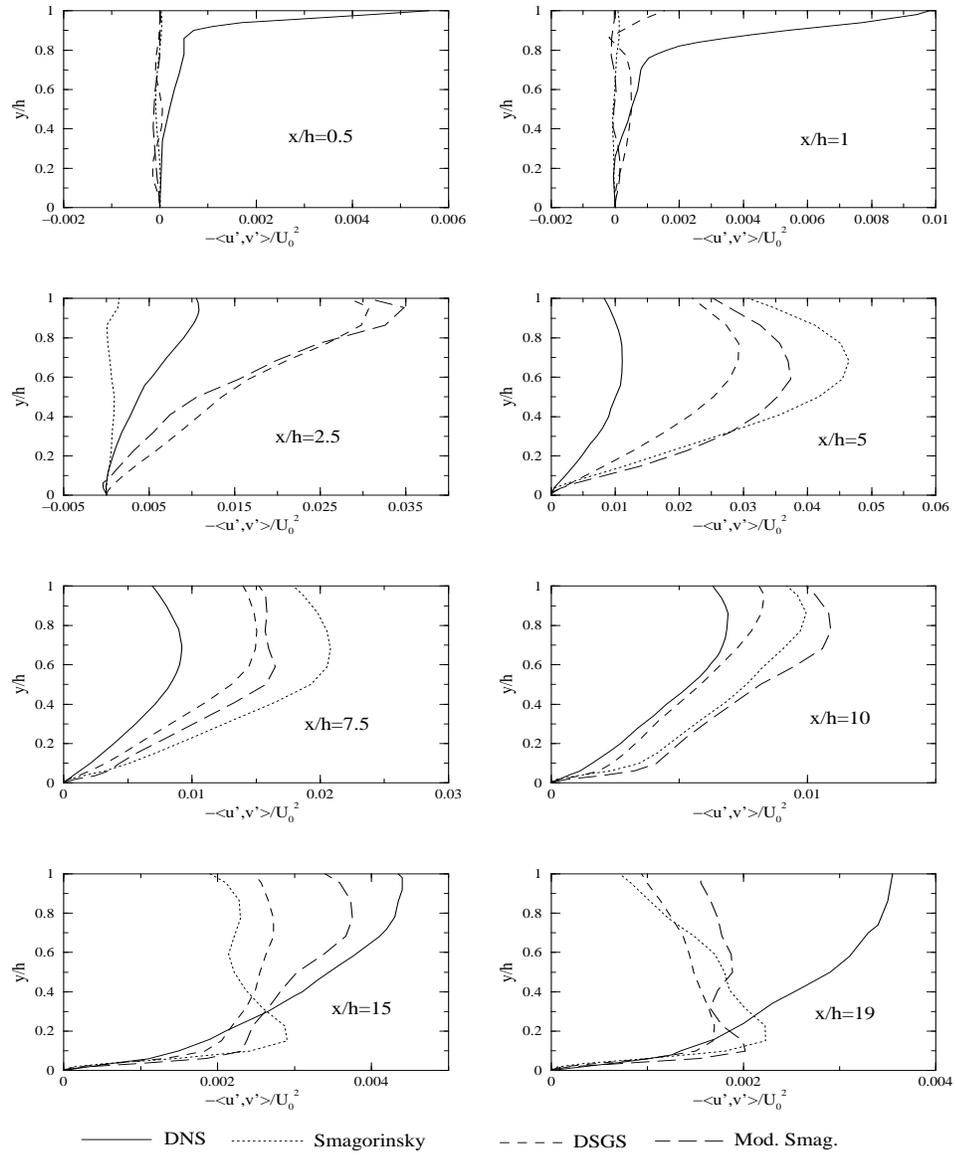}
\end{center}
\renewcommand{\baselinestretch}{1}
\caption{ \label{bfs-uv} \em
          Reynolds shear stress profiles in the post 
          expansion section of the backward facing step problem. DNS:
          direct numercial simulation data of Le and Moin \cite{le},
          Smagorinsky: original Smagorinsky subgrid scale model, DSGS: dynamic
          subgrid scale model, and Mod.~Smagorinsky: modified Smagorinsky model.}
\end{figure}

Figures \ref{bfs-u-rms}--\ref{bfs-uv} are the profiles of the various
Reynolds stress components. Figure \ref{bfs-u-rms}
indicates that in the reattachment region, all the SGS models predicted
a greater streamwise turbulence intensity near the wall than the DNS
results. According to Le and Moin \cite{le}, their
results were lower than the experimental data of Jovic and Driver \cite{jovic}.
Also, the data from Akselvoll and Moin \cite{akselvoll}, indicate that
their LES also predicted higher turbulence intensities, although not as 
high as our results. In the recirculation region, Le and Moin \cite{le}
reported that 
the maximum of all Reynolds stress components occurred around the step height. 
As can be seen in Figures \ref{bfs-u-rms}--\ref{bfs-uv}, the maximum of
all Reynolds stress components in the LES cases occur before the step
height and does not peak as sharply as in the DNS case.

\clearpage
\newpage

\section{Computer Software}

The numerical package used in this study is known as Hydra. This software
was developed by Mark Christon for Lawrence Livermore National 
Laboratory. The finite element code was in turn developed on the works
of Gresho et al.~\cite{gresho1}, \cite{gresho2}, \cite{gresho3}. 
The original Smagorinsky turbulence model was incorporated into the code
based on the work of Rose McCallen. Also a disscussion of the motivation 
for using the finite element method is also presented by McCallen
\cite{mccallen}. 

Although Hydra has several options for solving the
Pressure Poisson Equation, it was found that for our turbulent flow cases
the only solver that worked within a reasonable period of time was a
version of the direct solver. This together with our limited computer
resources restricted our mesh size to the $88\times 32 \times 16$ grid that
was used for the above results.

\clearpage
\newpage

\section{Conclusions}
One of the unique characteristics of LES as compared to other methods of 
dealing with turbulent flows is the idea of filtering. As was mentioned before,
this concept of filtering introduces a new error to the system which is the
result of the noncommutative nature of the filtering operation and 
differentiation in graded computational grids. Most of the focus in studying
this error has been on the filter kernal function $G$ and trying to obtain
a $G$ that would allow control over the order of the interchange error.
As we used FEM as our numerical method, we had knowledge of the general shape
of the solution (element wise) and so decided to focus on trying to use this
extra information. We were able to show that under certain conditions the
interchange error is second order. 

As we have noted previously, one of the weaknesses of the Smagorinsky 
subgrid scale model is its inability to damp out the eddy viscosity near
walls. The purpose of this project was to develope and implement a new
subgrid scale model that would address this problem. To test the new model
two turbulent flows where studied, the channel flow and the backward facing
step. The channel flow was chosen because of the extensive amount of data,
experimental, computational, and analytical, that is available for that 
type of flow. The backward facing step was chosen because of its geometric
simplicity for a flow with separation and reattachment.
Also, the Smagorinsky and a version of the 
Dynamic Subgird Scale model were also implemented for comparison.
These two models were chosen for comparison for two reasons. First, they
represent the two SGS models that are in common use in LES and secondly, their 
treatment of the eddy viscosity near walls are very different. As mentioned
above, the Smagorinsky model does not damp out the eddy viscosity near the
walls whereas the DSGS model does damp them out. Although the DSGS model
does damp out eddy viscosity near the walls, some problems with the model,
as mentioned previously, are that it requires double filtering and the
eddy viscosity becomes negative and leads to numerical instability if left
unchecked.

The results of our computations with the modified Smagorinsky model were
comparable to the DSGS model computations and generally performed better
than the original Smagorinsky model. Our eddy viscosity data from the 
channel flow indicate that our model and the DSGS were indeed able to 
detect the walls of the channel and reduce the eddy viscosity accordingly.

The next step in the continuation of this work should be to develope a
better numerical solver that would allow a finer grid and, in the case
of the backward facing step, allow an inlet channel to also be used.

\clearpage
\newpage
{\renewcommand{\baselinestretch}{1} \large \normalsize
\bibliographystyle{plain}
\bibliography{refs}
}
\newpage

\appendix

\section{Filter/Derivative Interchange Error}
\label{fltr-err-comp}

In this section we analyze the error in the interchange of filtering and
differentiation by numerically calculating the filter of the derivative
and the derivative of the filter and comparing the difference between the two.
We begin by first creating our graded mesh. The mesh is generated by taking 
an interval $[a,b]$, in our case $[0,1]$, and partioning it into n even 
intervals. The graded mesh is then generated by using the following
transformation:
\[
   x_i = \frac{e^{\alpha(a+ih)}-1}{e^\alpha -1}
\]
where $h=\frac{1}{n}$ and $\alpha$ is some constant, which in our case is
$\alpha=5$. If we define the filter width, $\del(x)$, to be twice the length
of the smallest element containing $x$, then we have
\begin{eqnarray*}
  \del(x_i) &=& 2(x_i - x_{i-1})
  \\
  &=& 2\left[\frac{e^{\alpha(a+ih)}-1}{e^\alpha-1}
             - \frac{e^{\alpha(a+(i-1)h)}-1}{e^\alpha-1}\right]
  \\
  &=& 2\left[\frac{e^{\alpha(a+ih)}(1-e^{-\alpha h})}{e^\alpha-1}\right]
  \\
  &=& 2\left[\frac{(e^{\alpha(a+ih)}-1+1)(1-e^{-\alpha h})}{e^\alpha-1}\right]
  \\
  &=& 2\left[x_i(1-e^{-\alpha h}) + \frac{1-e^{-\alpha h}}{e^\alpha-1}\right].
\end{eqnarray*}
So
\begin{equation}
\label{fltr-width}
   \del(x) = 2\left[x(1-e^{-\alpha h}) 
           + \frac{1-e^{-\alpha h}}{e^\alpha-1}\right]
\end{equation}
and
\begin{equation}
\label{deriv-fltr-width}
   \del'(x) = 2(1-e^{-\alpha h}).
\end{equation}
Taking the Taylor expansion of the exponential function about 0 in
Eqn.~(\ref{fltr-width}) and (\ref{deriv-fltr-width}), it is clear that the
order of the filter width and its derivative is $h$, i.~e. $O(\del(x))=
O(\del'(x))=O(h)$.

The numerical differentiation is performed using a fourth order method and 
the integration is carried out using the Composite Simpson's Rule over each
element in the filter support. This ensures that the numerical integration is
also at least fourth order in $h$. We use a truncated Gaussian as our filter
and two different filter widths were tested, $\del_1(x_i)= 2(x_i-x_{i-1})$ and
$\del_2(x_i)=2(x_i-x_{i-2})$.

To test the validity of our assertion that the numerical differentiation and
integration methods used are in fact fourth order, we observe the following.
In section \ref{fltr-fem-sec}, page \pageref{fem-deriv-err} we found the
interchange error to be
\[
   \left|\frac{d\ol{u}}{dx}(x)-\ol{\frac{du}{dx}}(x)\right| =
     \left|\del'(x)\int^{\frac{1}{2}}_{-\frac{1}{2}}
               yG(y)u'_e(x-y\del(x)) \,dy.\right|.
\]
Note that if our original function is a line, e.~g. $u(x)=10x+1$, then our 
polynomial interpolations will all be exact, i.~e. $u_e(x)=u(x)$ everywhere. 
Furthermore, the derivative is constant and the same in every element, e.~g.
for the above linear $u$, $u'_e(x) = 10$. So the interchange error should be
0. Hence the only error observed with the linear $u$ will be the error in the
calculation of the derivatives and integrals. 
\begin{table}[b]
\begin{center}
\begin{tabular}{|l||c|c||c|c|}\hline
\multicolumn{5}{|c|}{\rule[-3mm]{0mm}{8mm} 
 $\left\|\frac{d\ol{u}}{dx}(x)-\ol{\frac{du}{dx}}(x)\right\|_\infty$} \\
 \hline
 & \multicolumn{2}{|c||}{linear interpolation} 
 & \multicolumn{2}{|c|}{quadratic interpolation}
 \\ \hline
 h & $\del_1(x)$ & $\del_2(x)$ & $\del_1(x)$ & $\del_2(x)$
 \\ \hline\hline
 0.05 & 1.20e-03 & 1.20e-03 & 1.20e-03 & 1.20e-03
 \\ \hline
 0.005 & 1.19e-07 & 1.06e-07 & 1.19e-07 & 1.06e-07
 \\ \hline
 0.0005 & 1.72e-10 & 1.68e-10 & 1.70e-10 & 1.68e-10
 \\ \hline
\end{tabular}
\end{center}
\renewcommand{\baselinestretch}{1}
\caption{\label{fltr-err-table-lin}\em
    Max norm of the diff./filter interchange error, linear source function
    $u(x)=10x+1$}
\end{table}
Table~\ref{fltr-err-table-lin} indicates fourth order behavior and so our
claim of fourth order accuracy of the numerical differentiation and
integration is justified. 

We test the interchange error by using a sinusoidal source function
\[
   u(x)=\sin(2\pi x)+1.
\]
Table~\ref{fltr-err-table-sin} clearly indicates that
the interchage error is second order. 
\begin{table}
\begin{center}
\begin{tabular}{|l||c|c||c|c|}\hline
\multicolumn{5}{|c|}{\rule[-3mm]{0mm}{8mm} 
 $\left\|\frac{d\ol{u}}{dx}(x)-\ol{\frac{du}{dx}}(x)\right\|_\infty$} \\
 \hline
 & \multicolumn{2}{|c||}{linear interpolation} 
 & \multicolumn{2}{|c|}{quadratic interpolation}
 \\ \hline
 h & $\del_1(x)$ & $\del_2(x)$ & $\del_1(x)$ & $\del_2(x)$
 \\ \hline\hline
 0.01 & 1.07e-02 & 4.31e-02 & 1.37e-02 & 4.94e-02
 \\ \hline
 0.001 & 1.20e-04 & 4.97e-04 & 1.55e-04 & 6.13e-04
 \\ \hline
 0.0001 & 1.60e-06 & 5.00e-06 & 1.56e-06 & 6.18e-06
 \\ \hline
\end{tabular}
\end{center}
\renewcommand{\baselinestretch}{1}
\caption{\label{fltr-err-table-sin}\em
    Max norm of the diff./filter interchange error, sine source function
    $u(x)=\sin(2\pi x)+1$}
\end{table}
Hence our claim that the interchange
error is second order if $O(\del'(x))=O(\del(x))$ is justified.

\clearpage
\newpage

\section{The Effect of Double Filtering}
\label{double-filtering}

The following is a verification of the assumption, in Germano's derivation
of his Dynamic Subgrid Scale Model, that filtering twice is equivalent to
filtering once with some filter. In this paper we look at two filters that
are often mentioned with regard to LES. Namely the Gaussian filter and the
Top-Hat filter. 

\subsection{Gaussian Filter}

Given a Gaussian function
\[
  G(t) = e^{-\alpha t^2} = \frac{1}{2\pi}\int G(k) e^{ikt} \,dk
\]
the Fourier Transform of $G(t)$ denoted as $G(k)$ is
\[
  G(k) = \sqrt{\frac{\pi}{\alpha}} e^{\frac{-k^2}{4\alpha}}
       = \int G(t) e^{-ikt} \,dt
\]
Now we look at the following integral,
\[
  \ol{u}(x) = \int u(\xi) G(x-\xi) d\xi = u(x)*G(x)
\]
which is a more general case of the familiar Gaussian filter in LES. It 
will be easier to insert the parameters later to change this into the 
proper Gaussian filter. Let
\[
  G_1(t) = e^{-\alpha_1 t^2}
\]
and 
\[
  G_2(t) = e^{-\alpha_2 t^2}
\]
then we define
\[
  \ol{u}(x) = \int u(\xi) G_1(x-\xi) d\xi = u(x)*G_1(x)
\]
and
\[
  \wdtld{\ol{u}}(x) = \ol{u}(x)*G_2(x)
\]
So the Time Convolution Theorem, which states that the Fourier Transform
of the convolution of two functions is equal to the product of the Fourier
Transform of each of the functions, gives us
\[
  \ol{u}(k) = u(k) G_1(k)
\]
and
\begin{eqnarray*}
  \wdtld{\ol{u}}(k) &=& \ol{u}(k) G_2(k) 
\\
   &=& u(k) G_1(k) G_2(k)
\\
   &=& \frac{\pi}{\sqrt{\alpha_1 \alpha_2}} u(k) 
       e^{\frac{-k^2}{4}\left(\frac{1}{\alpha_1} + \frac{1}{\alpha_2}\right)}
\end{eqnarray*}
Hence using the inverse transform, we get
\begin{eqnarray*}
  \wdtld{\ol{u}}(x) &=& \frac{\pi}{\sqrt{\alpha_1 \alpha_2}} \frac{1}{2\pi}
     \int u(k) e^{\frac{-k^2}{4}\left(\frac{1}{\alpha_1} +
                                      \frac{1}{\alpha_2}\right)}
          e^{ikx} \, dk
\\
  &=& \frac{\pi}{\sqrt{\alpha_1 \alpha_2}} \frac{1}{2\pi}
      \int e^{\frac{-k^2}{4}\left(\frac{1}{\alpha_1} +
                                  \frac{1}{\alpha_2}\right)} e^{ikx}
      \int u(\xi) e^{-ik\xi} \,d\xi dk
\\
  &=& \frac{1}{2\sqrt{\alpha_1 \alpha_2}} \int u(\xi)
      \int e^{\frac{-k^2}{4}\left(\frac{1}{\alpha_1} + 
                                  \frac{1}{\alpha_2}\right)}
           e^{ik(x-\xi)} \, dk d\xi
\\
\end{eqnarray*}
Now if we let 
\[ \frac{1}{\alpha_3} = \frac{1}{\alpha_1} + \frac{1}{\alpha_2} \]
so
\[ \alpha_3 = \frac{\alpha_1 \alpha_2}{\alpha_1 + \alpha_2}\]
and
\[ 
  G_3(x) = e^{-\alpha_3 x^2} = \frac{1}{2\pi}\int G_3(k) e^{ikx} \,dk
\] 
then
\[
  G_3(k) = \sqrt{\frac{\pi}{\alpha_3}} e^{-\frac{k^2}{4\alpha_3}}
\]
So
\begin{eqnarray*}
  \wdtld{\ol{u}}(x) &=& \frac{1}{2\sqrt{\alpha_1 \alpha_2}} \int u(\xi)
                        \int e^{\frac{-k^2}{4}\left(\frac{1}{\alpha_1} 
                      + \frac{1}{\alpha_2}\right)} e^{ik(x-\xi)} \, dk d\xi
\\
  &=& \frac{1}{2 \sqrt{\pi}\sqrt{\alpha_1 + \alpha_2}}
    \int u(\xi) \int G_3(k) e^{ik(x-\xi)} \, dk d\xi
\\
  &=& \sqrt{\frac{\pi}{\alpha_1 + \alpha_2}} \int u(\xi) G_3(x-\xi) \,d\xi
\\
  &=& \sqrt{\frac{\pi}{\alpha_1 + \alpha_2}} 
      \int u(\xi) e^{-(x-\xi)^2 \alpha} \,d\xi
\\
\end{eqnarray*}

Now the Gaussian filter for LES is 
\[ G(x) = \sqrt{\frac{6}{\pi}} \frac{1}{\del} e^{-6 \frac{x^2}{\del^2}}\]
So if we let
\[ G_1(x) = \sqrt{\frac{6}{\pi}} \frac{1}{\del_1} e^{-6 \frac{x^2}{\del_1^2}}\]
\[ G_2(x) = \sqrt{\frac{6}{\pi}} \frac{1}{\del_2} e^{-6 \frac{x^2}{\del_2^2}}\]
then $\alpha_1 = \frac{6}{\del_1^2}$, $\alpha_2 = \frac{6}{\del_2^2}$ and
$\alpha_3 = \frac{6}{\del_1^2 + \del_2^2}$.
Hence if we go back and let $\ol{\cdot}$ and $\wdtld{\cdot}$ be the proper
Gaussian filters defined with the above kernals, we have that 
\begin{eqnarray*}
 \wdtld{\ol{u}}(x) &=& [u(x)*G_1(x)]*G_2(x)
\\
   &=& \frac{6}{\pi} \frac{\sqrt{\pi}}{\del_1 \del_2} 
       \frac{\del_1 \del_2}{\sqrt{6(\del^2_1 + \del^2_2}}
       \int u(\xi) e^{-\frac{(x-\xi)^2}{\del^2_1 + \del^2_2}} \,d\xi
\\
   &=& \sqrt{\frac{6}{\pi}} \frac{1}{\sqrt{\del_1^2 + \del_2^2}}
       \int u(\xi) e^{-6 \frac{(x-\xi)^2}{\del_1^2 + \del_2^2}} \,d\xi
\\
  &=& u(x)*G_s(x)
\end{eqnarray*}
where $G_s(x)$ is the Gaussian filter with filter width 
$\del = \sqrt{\del_1^2 + \del_2^2}$.

\subsection{Top--Hat Filter}

The Top--Hat filter is defined as 
\[
   p_{\frac{\del}{2}}(x) = \left\{ \begin{array}{r@{\quad:\quad}l}
                          \frac{1}{\del} & |x| < \frac{\del}{2}
                      \\  0 & \mbox{otherwise}
                          \end{array} \right. 
\]
We also define a triangle filter as follows.
\[
  q_{\del}(x) = \frac{1}{\del}
                \left\{ \begin{array}{r@{\quad:\quad}l}
                1 + \frac{x}{\del} & -\del < x < 0
            \\  1 - \frac{x}{\del} & 0 < x < \del
            \\  0 & \mbox{otherwise}
              \end{array} \right.
\]
Now the Fourier Transform of the above functions are
\[
  p_{\frac{\del}{2}}(k) = \frac{2 \sin(\frac{k \del}{2})}{k \del}
\]
and
\[
  q_\del(k) = \frac{4 \sin^2(\frac{k \del}{2})}{\del^2 k^2}
\]
If we define $\ol{\cdot}$ to be the Top-Hat filter and $\wdtld{\cdot}$ to be
the triangle filter, then note that if we filter twice with the same (width)
Top-Hat filter we get
\[
  \ol{\ol{u}}(x) = u * p_{\frac{\del}{2}} * p_{\frac{\del}{2}}
\]
So the Fourier Transform of $\ol{\ol{u}}$ is
\begin{eqnarray*}
  \ol{\ol{u}}(k) &=& u(k) p_{\frac{\del}{2}}(k) p_{\frac{\del}{2}}(k)
\\
                 &=& u(k) \frac{4 \sin^2(\frac{k \del}{2})}{k^2 \del^2}
\end{eqnarray*}
Hence 
\begin{eqnarray*}
  \ol{\ol{u}}(x) &=& \frac{1}{2\pi} \int \ol{\ol{u}}(k) e^{ikx} \, dk
\\
     &=& \frac{1}{2\pi \del^2} 
         \int \frac{4 \sin^2(\frac{k\del}{2})}{k^2} e^{ikx}
         \int u(\xi) e^{-ik\xi} \,d\xi dk
\\
     &=& \frac{1}{2\pi \del^2} \int u(\xi)
         \int \frac{4 \sin^2(\frac{k\del}{2})}{k^2} e^{ik(x-\xi)} \, dk d\xi
\\
     &=& \frac{1}{\del} \int u(\xi) \frac{1}{2\pi}
         \int \frac{4 \sin^2(\frac{k\del}{2})}{\del^2 k^2} e^{ik(x-\xi)} 
         \, dk d\xi
\\
     &=& \frac{1}{\del} \int u(\xi) q_\del(x-\xi) \, d\xi
\\
     &=& \wdtld{u}(x)
\end{eqnarray*}
where
\[
  q_\del(x-\xi) = \frac{1}{\del}
                  \left\{ \begin{array}{r@{\quad:\quad}l}
                  1 + \frac{\xi - x}{\del} & x-\del < \xi < x
              \\  1 - \frac{\xi - x}{\del} & x < \xi < x + \del
              \\  0 & \mbox{otherwise}
                  \end{array} \right.
\]

So filtering with the same Top--Hat filter twice is equivalent to 
filtering once with the triangle filter, $q_\del$.

\end{document}